\documentclass[12pt]{article}

 \usepackage{mathrsfs,amsfonts,amssymb}

 \newtheorem{lemma}{Lemma}[section]
 \newtheorem{theorem}{Theorem}[section]
 \newtheorem{corollary}{Corollary}[section]

 \def\blemma{\begin{lemma}\sl{}\def\elemma{\end{lemma}}}
 \def\btheorem{\begin{theorem}\sl{}\def\etheorem{\end{theorem}}}

 \def\beqlb{\begin{eqnarray}}\def\eeqlb{\end{eqnarray}}
 \def\beqnn{\begin{eqnarray*}}\def\eeqnn{\end{eqnarray*}}

 \def\qqquad{\qquad\qquad}
 \def\qed{\hfill{$\square$}\medskip}

 \def\<{\langle}\def\>{\rangle}

 \def\E{\mbox{\boldmath $E$}}
 \def\P{\mbox{\boldmath $P$}}\def\Q{\mbox{\boldmath $Q$}}

 \def\bB{\mbox{\boldmath $B$}}\def\bM{\mbox{\boldmath $M$}}

 \def\A{{\cal A}}\def\B{{\cal B}}\def\L{{\cal L}}

 \def\itDelta{{\it\Delta}}\def\itGamma{{\it\Gamma}}
 \def\itOmega{{\it\Omega}}\def\itPhi{{\it\Phi}}
 \def\itPsi{{\it\Psi}}

 \def\IR{{I\!\!R}}

\begin{document}

\noindent{Published in: {\it Electronic Journal of Probability}
{\bf 6} (2001), Paper No. 25, 1--33.}

\bigskip\bigskip

\centerline{\Large\bf Superprocesses with Dependent Spatial}

\medskip
\centerline{\Large\bf Motion and General Branching Densities}

\bigskip
\centerline{Donald A. Dawson\footnote{Supported by NSERC operating grant.}}

\smallskip
\centerline{School of Mathematics and Statistics, Carleton University,}

\centerline{1125 Colonel By Drive, Ottawa, Canada K1S 5B6}

\bigskip
\centerline{Zenghu Li\footnote{Supported by NNSF grant 19361060.}}

\smallskip
\centerline{Department of Mathematics, Beijing Normal University,}

\centerline{Beijing 100875, P.R. China}

\bigskip
\centerline{Hao Wang\footnote{Supported by the research grant of UO.}}

\smallskip
\centerline{Department of Mathematics, University of Oregon,}

\centerline{Eugene OR 97403-1222, U.S.A.}

\bigskip\bigskip

{\narrower{\narrower

\centerline{\bf Abstract}

\bigskip

We construct a class of superprocesses by taking the high density limit of a sequence
of interacting-branching particle systems. The spatial motion of the superprocess is
determined by a system of interacting diffusions, the branching density is given by an
arbitrary bounded non-negative Borel function, and the superprocess is characterized
by a martingale problem as a diffusion process with state space $M(\IR)$, improving
and extending considerably the construction of Wang (1997, 1998). It is then proved
in a special case that a suitable rescaled process of the superprocess converges to
the usual super Brownian motion. An extension to measure-valued branching catalysts is
also discussed.

\bigskip\bigskip


\bigskip


\par}\par}

\bigskip\bigskip


\section{Introduction}
\setcounter{equation}{0}

For a given topological space $E$, let $B(E)$ denote the totality of all bounded
Borel functions on $E$ and let $C(E)$ denote its subset comprising of continuous
functions. Let $M(E)$ denote the space of finite Borel measures on $E$ endowed with
the topology of weak convergence. Write $\<f,\mu\>$ for $\int fd\mu$. For $F\in
B(M(E))$ let
 \beqlb\label{1.1}
\frac{\delta F(\mu)}{\delta\mu(x)}
=
\lim_{r\to 0^+}\frac{1}{\,r\,}[F(\mu + r\delta_x) - F(\mu)],
\qquad x\in E,
 \eeqlb
if the limit exists. Let $\delta^2F(\mu) / \delta\mu(x)\delta\mu(y)$ be defined in
the same way with $F$ replaced by $(\delta F/ \delta\mu(y))$ on the right hand side.
For example, if $F_{m,f}(\mu) = \<f,\mu^m\>$ for $f\in B(E^m)$ and $\mu \in M(E)$,
then
 \beqlb\label{1.2}
\frac{\delta F_{m,f}(\mu)}{\delta\mu(x)}
=
\sum_{i=1}^m\<\itPsi_i (x)f,\mu^{m-1}\>,
\qquad x\in E,
 \eeqlb
where $\itPsi_i(x)$ is the operator from $B(E^m)$ to $B(E^{m-1})$ defined by
 \beqlb\label{1.3}
\itPsi_i(x)f(x_1,\cdots,x_{m-1}) = f(x_1,\cdots,x_{i-1},x,x_i,\cdots,x_{m-1}),
\quad x_j\in E,
 \eeqlb
where $x\in E$ is the $i$th variable of $f$ on the right hand side.

Now we consider the case where $E=\IR$, the one-dimensional Euclidean space. Suppose
that $c\in C(\IR)$ is Lipschitz and $h\in C(\IR)$ is square-integrable. Let
 \beqlb\label{1.4}
\rho(x) = \int_{\IR}h(y-x)h(y)dy,
 \eeqlb
and $a(x) = c(x)^2 + \rho(0)$ for $x\in \IR$. We assume in addition that $\rho$ is
twice continuously differentiable with $\rho^{\prime}$ and $\rho^{\prime\prime}$
bounded, which is satisfied if $h$ is integrable and twice continuously
differentiable with $h^\prime$ and $h^{\prime\prime}$ bounded. Then
 \beqlb\label{1.5}
\A F(\mu)
&=&
\frac{1}{2}\int_{\IR}a(x)\frac{d^2}{dx^2}
\frac{\delta F(\mu)}{\delta\mu(x)}\mu(dx)   \nonumber  \\
& &
+\,\frac{1}{2}\int_{\IR^2}\rho(x-y)\frac{d^2}{dxdy}
\frac{\delta^2 F(\mu)}{\delta\mu(x)\delta\mu(y)}\mu(dx)\mu(dy)
 \eeqlb
defines an operator $\A$ which acts on a subset of $B(M(\IR))$ and generates a
diffusion process with state space $M(\IR)$. Suppose that $\{W(x,t): x\in \IR,
t\ge 0\}$ is a Brownian sheet and $\{B_i(t): t\ge0\}$, $i=1,2,\cdots$, is a family of
independent standard Brownian motions which are independent of $\{W(x,t): x\in \IR,
t\ge 0\}$. By Lemma \ref{l3.1}, for any initial conditions $x_i(0) = x_i$, the
stochastic equations
 \beqlb\label{1.6}
dx_i(t)
= c(x_i(t)) dB_i(t)
+ \int_{\IR} h(y-x_i(t))W(dy, dt),
\qquad t\ge0, i=1,2,\cdots,
 \eeqlb
have unique solutions $\{x_i(t): t\ge 0\}$ and, for each integer $m\ge1$, $\{(x_1(t),
\cdots, x_m(t)): t\ge0\}$ is an $m$-dimensional diffusion process which is generated
by the differential operator
 \beqlb\label{1.7}
G^m :=
\displaystyle{
\frac{1}{2}\sum_{i=1}^m
a(x_i)\frac{\partial^{2}}{\partial x_{i}^{2}}}
\displaystyle{ + \frac{1}{2} \sum_{i,j=1, i\neq j}^{m}
\rho(x_{i}-x_{j})\frac{\partial^{2}}{\partial x_{i} \partial x_{j}}}.
 \eeqlb
In particular, $\{x_i(t): t\ge0\}$ is a one-dimensional diffusion process with
generator $G:= (a(x)/2) \itDelta$. Because of the exchangebility, a diffusion process
generated by $G^m$ can be regarded as an interacting particle system or a
measure-valued process. Heuristically, $a(\cdot)$ represents the speed of the
particles and $\rho(\cdot)$ describes the interaction between them. The diffusion
process generated by $\A$ arises as the high density limit of a sequence of
interacting particle systems described by (\ref{1.6}); see Wang (1997, 1998) and
section 4 of this paper. For $\sigma \in B(\IR)^+$, we may also define the operator
$\B$ by
 \beqlb\label{1.8}
\B F(\mu)
=
\frac{1}{2}\int_{\IR}\sigma(x)
\frac{\delta^2 F(\mu)}{\delta\mu(x)^2}\mu(dx).
 \eeqlb
A Markov process generated by $\L := \A + \B$ is naturally called a {\it superprocess
with dependent spatial motion (SDSM)} with parameters $(a, \rho, \sigma)$, where
$\sigma$ represents the branching density of the process. In the special case where
both $c$ and $\sigma$ are constants, the SDSM was constructed in Wang (1997, 1998) as a
diffusion process in $M(\hat\IR)$, where $\hat \IR = \IR \cup \{\partial\}$ is the
one-point compactification of $\IR$. It was also assumed in Wang (1997, 1998) that $h$
is a symmetric function and that the initial state of the SDSM has compact support in
$\IR$. Stochastic partial differential equations and local times associated with the
SDSM were studied in Dawson et al (2000a, b).

The SDSM contains as special cases several models arising in different circumstances
such as the one-dimensional super Brownian motion, the molecular diffusion with
turbulent transport and some interacting diffusion systems of McKean-Vlasov type;
see e.g. Chow (1976), Dawson (1994), Dawson and Vaillancourt (1995) and Kotelenez
(1992, 1995). It is thus of interest to construct the SDSM under reasonably more
general conditions and formulate it as a diffusion processes in $M(\IR)$. This is the
main purpose of the present paper. The rest of this paragraph describes the main
results of the paper and gives some unsolved problems in the subject. In section 2, we
define some function-valued dual process and investigate its connection to the
solution of the martingale problem of a SDSM. Duality method plays an important role in
the investigation. Although the SDSM could arise as high density limit of a sequence
of interacting-branching particle systems with location-dependent killing density
$\sigma$ and binary branching distribution, the construction of such systems seems
rather sophisticated and is thus avoided in this work. In section 3, we construct the
interacting-branching particle system with uniform killing density and
location-dependent branching distribution, which is comparatively easier to treat. The
arguments are similar to
those in Wang (1998). The high density limit of the interacting-branching particle
system is considered in section 4, which gives a solution of the martingale problem of
the SDSM in the special case where $\sigma \in C(\IR)^+$ can be extended into a
continuous function on $\hat\IR$. In section 5, we use the dual process to extend the
construction of the SDSM to a general bounded Borel branching density $\sigma \in
B(\IR)^+$. In both sections 4 and 5, we use martingale arguments to show that, if the
processes are initially supported by $\IR$, they always stay in $M(\IR)$, which are
new results even in the special case considered in Wang (1997, 1998). In section 6, we
prove a rescaled limit theorem of the SDSM, which states that a suitable rescaled
SDSM converges to the usual super Brownian motion if $c(\cdot)$ is bounded away from
zero. This describes another situation where the super Brownian motion arises
universally; see also Durrett and Perkins (1998) and Hara and Slade (2000a, b). When
$c(\cdot) \equiv 0$, we expect that the same rescaled limit would lead to a
measure-valued diffusion process which is the high density limit of a sequence of
coalescing-branching particle systems, but there is still a long way to reach a
rigorous proof. It suffices to mention that not only the characterization of those
high density limits but also that of the coalescing-branching particle systems
themselves are still open problems. We refer the reader to Evans and Pitman (1998)
and the references therein for some recent work on related models. In section 7, we
consider an extension of the construction of the SDSM to the case where $\sigma$ is
of the form $\sigma = \dot\eta$ with $\eta$ belonging to a large class of Radon
measures on $\IR$, in the lines of Dawson and Fleischmann (1991, 1992). The process
is constructed only when $c(\cdot)$ is bounded away from zero and it can be called a
{\it SDSM with measure-valued catalysts}. The transition semigroup of the SDSM with
measure-valued catalysts is constructed and characterized using a measure-valued dual
process. The derivation is based on some estimates of moments of the dual process.
However, the existence of a diffusion realization of the SDSM with measure-valued
catalysts is left as another open problem in the subject.

Notation: Recall that $\hat \IR = \IR \cup \{\partial\}$ denotes the one-point
compactification of $\IR$. Let $\lambda^m$ denote the Lebesgue measure on $\IR^m$. Let
$C^2(\IR^m)$ be the set of twice continuously differentiable functions on
$\IR^m$ and let $C^2_\partial (\IR^m)$ be the set of functions in $C^2 (\IR^m)$
which together with their derivatives up to the second order can be extended
continuously to $\hat \IR$. Let $C_0^2 (\IR^m)$ be the subset of $C_\partial^2
(\IR^m)$ of functions that together with their derivatives up to the second order {\it
vanish rapidly} at infinity. Let $(T_t^m)_{t\ge0}$ denote the transition semigroup of
the $m$-dimensional standard Brownian motion and let $(P_t^m)_{t\ge0}$ denote the
transition semigroup generated by the operator $G^m$. We shall omit the superscript
$m$ when it is one. Let $(\hat P_t)_{t\ge0}$ and $\hat G$ denote the extensions of
$(P_t)_{t\ge0}$ and $G$ to $\hat \IR$ with $\partial$ as a trap. We denote the
expectation by the letter of the probability measure if this is specified and simply
by $\E$ if the measure is not specified.

We remark that, if $|c(x)| \ge \epsilon >0$ for all $x\in \IR$, the semigroup
$(P_t^m)_{t>0}$ has density $p_t^m(x,y)$ which satisfies
 \beqlb\label{1.9}
p_t^m(x,y) \le \mbox{const}\cdot g_{\epsilon t}^m(x,y),
\qquad t>0, x,y\in \IR^m,
 \eeqlb
where $g_t^m(x,y)$ denotes the transition density of the $m$-dimensional standard
Brownian motion; see e.g. Friedman (1964, p.24).


\section{Function-valued dual processes}
\setcounter{equation}{0}
In this section, we define a function-valued dual process and investigate its
connection to the solution of the martingale problem for the SDSM. Recall the
definition of the generator $\L := \A +\B$ given by (\ref{1.5}) and (\ref{1.8}) with
$\sigma \in B(\IR)^+$. For $\mu\in M(\IR)$ and a subset ${\cal D}({\cal L})$ of the
domain of ${\cal L}$, we say an $M(\IR)$-valued c\'adl\'ag process $\{X_t: t\ge0\}$ is
a solution of the {\it $({\cal L}, {\cal D}({\cal L}),\mu)$-martingale problem} if
$X_0=\mu$ and
 \beqnn
F(X_t) - F(X_0) - \int_0^t {\cal L}F(X_s) ds,
\qquad t\ge0,
 \eeqnn
is a martingale for each $F\in{\cal D}({\cal L})$. Observe that, if $F_{m,f}(\mu) =
\<f,\mu^m\>$ for $f\in C^2(\IR^m)$, then
 \beqlb\label{2.1}
\A F_{m,f}(\mu)
&=&
\frac{1}{2}\int_{\IR^m} \sum_{i=1}^m a(x_i)
f^{\prime\prime}_{ii} (x_1,\cdots,x_m) \mu^m(dx_1,\cdots,dx_m)   \nonumber  \\
& &\quad + \frac{1}{2}\int_{\IR^m} \sum_{i,j=1, i\neq j}^m \rho(x_i-x_j)
f^{\prime\prime}_{ij} (x_1,\cdots,x_m) \mu^m(dx_1,\cdots,dx_m)  \nonumber  \\
&=&
F_{m,G^mf}(\mu),
 \eeqlb
and
 \beqlb\label{2.2}
\B F_{m,f}(\mu)
&=&
\frac{1}{2} \sum_{i,j=1, i\neq j}^m\int_{\IR^{m-1}} \itPhi_{ij}
f(x_1,\cdots,x_{m-1}) \mu^{m-1}(dx_1,\cdots,dx_{m-1})   \nonumber  \\
&=&
\frac{1}{2} \sum_{i,j=1, i\neq j}^m F_{m-1,\itPhi_{ij}f}(\mu),
 \eeqlb
where $\itPhi_{ij}$ denotes the operator from $B(\IR^m)$ to $B(\IR^{m-1})$
defined by
 \beqlb\label{2.3}
\itPhi_{ij}f(x_1,\cdots,x_{m-1})
=
\sigma(x_{m-1})f(x_1,\cdots,x_{m-1},\cdots,x_{m-1},\cdots,x_{m-2}),
 \eeqlb
where $x_{m-1}$ is in the places of the $i$th and the $j$th variables of $f$ on
the right hand side. It follows that
 \beqlb\label{2.4}
\L F_{m,f}(\mu)
=
F_{m,G^mf}(\mu)
+ \frac{1}{2} \sum_{i,j=1, i\neq j}^m F_{m-1,\itPhi_{ij}f}(\mu).
 \eeqlb

Let $\{M_t: t\ge 0\}$ be a nonnegative integer-valued c\'adl\'ag Markov process with
transition intensities $\{q_{i,j}\}$ such that $q_{i,i-1} = -q_{i,i} = i(i-1)/2$ and
$q_{i,j}=0$ for all other pairs $(i,j)$. That is, $\{M_t: t\ge 0\}$ is the well-known
Kingman's coalescent process. Let $\tau_0 =0$ and $\tau_{M_0} = \infty$, and let
$\{\tau_k: 1\le k\le M_0-1\}$ be the sequence of jump times of $\{M_t: t\ge 0\}$. Let
$\{\itGamma_k: 1\le k\le M_0-1\}$ be a sequence of random operators which are
conditionally independent given $\{M_t: t\ge 0\}$ and satisfy
 \beqlb\label{2.5}
\P\{\itGamma_k = \itPhi_{i,j} | M(\tau_k^-) =l\}
=
\frac{1}{l(l-1)},
\qquad 1 \le i \neq j \le l,
 \eeqlb
where $\itPhi_{i,j}$ is defined by (\ref{2.3}). Let $\bB$ denote the topological union
of $\{B(\IR^m): m=1,2, \cdots\}$ endowed with pointwise convergence on each
$B(\IR^m)$. Then
 \beqlb\label{2.6}
Y_t
=
P^{M_{\tau_k}}_{t-\tau_k} \itGamma_k P^{M_{\tau_{k-1}}}_{\tau_k -\tau_{k-1}}
\itGamma_{k-1} \cdots P^{M_{\tau_1}}_{\tau_2 -\tau_1} \itGamma_1
P^{M_0}_{\tau_1}Y_0,
\quad \tau_k \le t < \tau_{k+1}, 0\le k\le M_0-1,
 \eeqlb
defines a Markov process $\{Y_t: t\ge0\}$ taking values from $\bB$. Clearly, $\{(M_t,
Y_t): t\ge 0\}$ is also a Markov process. To simplify the presentation, we shall
suppress the dependence of $\{Y_t: t\ge 0\}$ on $\sigma$ and let $\E^\sigma_{m,f}$
denote the expectation given $M_0=m$ and $Y_0=f \in C(\IR^m)$, just as we are working
with a canonical realization of $\{(M_t, Y_t): t\ge 0\}$. By (\ref{2.6}) we have
 \beqlb\label{2.7}
& &\E^{\sigma}_{m,f} \bigg[\<Y_t, \mu^{M_t}\>
\exp\bigg\{\frac{1}{2}\int_0^t M_s(M_s-1)ds \bigg\}\bigg]  \nonumber \\
&=&
\<P^m_tf,\mu^m\>  \\
& &
+ \frac{1}{2}\sum_{i,j=1, i\neq j}^m \int_0^t
\E^{\sigma}_{m-1,\itPhi_{ij}P^m_uf} \bigg[\<Y_{t-u},\mu^{M_{t-u}}\>
\exp\bigg\{\frac{1}{2}\int_0^{t-u} M_s(M_s-1)ds \bigg\}\bigg]du.  \nonumber
 \eeqlb

\blemma\label{l2.1}
For any $f \in B(\IR^m)$ and any integer $m\ge1$,
 \beqlb\label{2.8}
& &\E^{\sigma}_{m,f} \bigg[\<Y_t, \mu^{M_t}\>
\exp\bigg\{\frac{1}{2}\int_0^t M_s(M_s-1)ds \bigg\}\bigg]   \nonumber \\
&\le&
\|f\| \sum_{k=0}^{m-1} 2^{-k}m^k(m-1)^k \|\sigma\|^k\<1,\mu\>^{m-k},
 \eeqlb
where $\|\cdot\|$ denotes the supremum norm.
\elemma

{\it Proof.} The left hand side of (\ref{2.8}) can be decomposed as $\sum_{k=0}^{m-1}
A_k$ with
 \beqnn
A_k = \E^{\sigma}_{m,f} \bigg[\<Y_t, \mu^{M_t}\>
\exp\bigg\{\frac{1}{2}\int_0^t M_s(M_s-1)ds \bigg\}
1_{\{\tau_k\le t<\tau_{k+1}\}}\bigg].
 \eeqnn
Observe that $A_0 = \<P^m_tf,\mu^m\> \le \|f\| \<1,\mu\>^m$ and
 \beqnn
A_k
&=&
\frac{m!(m-1)!}{2^k (m-k)!(m-k-1)!} \int_0^t ds_1 \int_{s_1}^t ds_2 \cdots   \\
& &\cdot\int_{s_{k-1}}^t \E^{\sigma}_{m,f}\{\<P^{m-k}_{t-s_k} \itGamma_k \cdots
P^{m-1}_{s_2 -s_1} \itGamma_1 P^m_{s_1}f,\mu^{m-k}\>
| \tau_j = s_j: 1\le j\le k\} ds_k   \\
&\le&
\frac{m!(m-1)!}{2^k (m-k)!(m-k-1)!} \int_0^t ds_1 \int_0^t ds_2 \cdots
\int_0^t \|f\|\|\sigma\|^k\<1,\mu\>^{m-k} ds_k   \\
&\le&
\frac{m!(m-1)!}{2^k (m-k)!(m-k-1)!}
\|f\|\|\sigma\|^k\<1,\mu\>^{m-k} t^k
 \eeqnn
for $1\le k\le m-1$. Then we get the conclusion. \qed

\blemma\label{l2.2}
Suppose that $\sigma_n \to \sigma$ boundedly and pointwise and $\mu_n \to \mu$ in
$M(\IR)$ as $n\to \infty$. Then, for any $f \in B(\IR^m)$ and any integer $m\ge1$,
 \beqlb\label{2.9}
& &\E^\sigma_{m,f} \bigg[\<Y_t, \mu^{M_t}\>
\exp\bigg\{\frac{1}{2}\int_0^t M_s(M_s-1)ds \bigg\}\bigg]  \nonumber   \\
&=&
\lim_{n\to \infty} \E^{\sigma_n}_{m,f} \bigg[\<Y_t, \mu_n^{M_t}\>
\exp\bigg\{\frac{1}{2}\int_0^t M_s(M_s-1)ds \bigg\}\bigg].
 \eeqlb
\elemma

{\it Proof.} For $h\in C(\IR^2)$ we see by (\ref{2.7}) that
 \beqlb\label{2.10}
& &\E^{\sigma_n}_{1,\itPhi_{12}h}\bigg[\<Y_t,\mu_n^{M_t}\>
\exp\bigg\{\frac{1}{2}\int_0^t M_s(M_s-1)ds \bigg\}\bigg]  \nonumber \\
&=&
\E^{\sigma_n}_{1,\itPhi_{21}h}\bigg[\<Y_t,\mu_n^{M_t}\>
\exp\bigg\{\frac{1}{2}\int_0^t M_s(M_s-1)ds \bigg\}\bigg]  \nonumber \\
&=&
\int_{\IR^2} h(y,y)p_t(x,y) \mu_n(dx)\sigma_n(y)dy.
 \eeqlb
If $f,g \in C(\IR)^+$ have bounded supports, then we have $f(x)\mu_n(dx) \to
f(x)\mu(dx)$ and $g(y) \sigma_n(y) dy \to g(y) \sigma(y) dy$ by weak convergence,
so that
 \beqnn
\lim_{n\to \infty}\int_{\IR^2} f(x)g(y)p_t(x,y) \mu_n(dx)\sigma_n(y)dy
=
\int_{\IR^2} f(x)g(y)p_t(x,y) \mu(dx)\sigma(y)dy.
 \eeqnn
Since $\{\mu_n\}$ is tight and $\{\sigma_n\}$ is bounded, one can easily see that
$\{p_t(x,y) \mu_n(dx)\sigma_n(y)dy\}$ is a tight sequence and hence $p_t(x,y)
\mu_n(dx) \sigma_n(y) dy \to p_t(x,y) \mu(dx) \sigma(y) dy$ by weak convergence.
Therefore, the value of (\ref{2.10}) converges as $n\to \infty$ to
 \beqnn
& &\E^{\sigma}_{1,\itPhi_{12}h}\bigg[\<Y_t,\mu^{M_t}\>
\exp\bigg\{\frac{1}{2}\int_0^t M_s(M_s-1)ds \bigg\}\bigg]  \nonumber \\
&=&
\E^{\sigma}_{1,\itPhi_{21}h}\bigg[\<Y_t,\mu^{M_t}\>
\exp\bigg\{\frac{1}{2}\int_0^t M_s(M_s-1)ds \bigg\}\bigg]  \\
&=&
\int_{\IR^2} h(y,y)p_t(x,y) \mu(dx)\sigma(y)dy.   \nonumber
 \eeqnn
Applying bounded convergence theorem to (\ref{2.7}) we get inductively
 \beqnn
& &\E^\sigma_{m-1,\itPhi_{ij}P^m_tf} \bigg[\<Y_t, \mu^{M_t}\>
\exp\bigg\{\frac{1}{2}\int_0^t M_s(M_s-1)ds \bigg\}\bigg]   \\
&=&
\lim_{n\to \infty}\E^{\sigma_n}_{m-1,\itPhi_{ij}P^m_tf}
\bigg[\<Y_t, \mu_n^{M_t}\>
\exp\bigg\{\frac{1}{2}\int_0^t M_s(M_s-1)ds \bigg\}\bigg]
 \eeqnn
for $1\le i\neq j\le m$. Then the result follows from (\ref{2.7}). \qed

\btheorem\label{t2.1}
Let ${\cal D}(\L)$ be the set of all functions of the form $F_{m,f}(\mu) = \<f,\mu^m\>$
with $f\in C^2 (\IR^m)$. Suppose that $\{X_t: t\ge0\}$ is a continuous
$M(\IR)$-valued process and that $\E \{\<1,X_t\>^m\}$ is locally bounded in $t\ge0$ for
each $m\ge1$. If $\{X_t: t\ge0\}$ is a solution of the $({\cal L}, {\cal D}(\L),
\mu)$-martingale problem, then
 \beqlb\label{2.11}
\E \<f,X^m_t\>
=
\E^\sigma_{m,f} \bigg[\<Y_t, \mu^{M_t}\>
\exp\bigg\{\frac{1}{2}\int_0^t M_s(M_s-1)ds \bigg\}\bigg]
 \eeqlb
for any $t\ge0$, $f \in B(\IR^m)$ and integer $m\ge1$.
\etheorem

{\it Proof.} In view of (\ref{2.6}), the general equality follows by bounded
pointwise approximation once it is proved for $f\in C^2 (\IR^m)$. In this
proof, we set $F_\mu(m,f) = F_{m,f}(\mu) = \<f,\mu^m\>$. From the construction
(\ref{2.6}), it is not hard to see that $\{(M_t, Y_t): t\ge 0\}$ has generator $\L^*$
given by
 \beqnn
\L^* F_\mu(m,f)
=
F_\mu(m,G^mf)
+ \frac{1}{2} \sum_{i,j=1, i\neq j}^m
[F_\mu(m-1,\itPhi_{ij}f) - F_\mu(m,f)].
 \eeqnn
In view of (\ref{2.4}) we have
 \beqlb\label{2.12}
\L^* F_\mu(m,f)
=
\L F_{m,f}(\mu) - \frac{1}{2} m(m-1) F_{m,f}(\mu).
 \eeqlb
The following calculations are guided by the relation (\ref{2.12}). In the sequel, we
assume that $\{X_t: t\ge 0\}$ and $\{(M_t, Y_t): t\ge 0\}$ are defined on the same
probability space and are independent of each other. Suppose that for each $n\ge1$ we
have a partition $\itDelta_n := \{0= t_0 <t_1 <\cdots <t_n =t\}$ of $[0,t]$. Let
$\|\itDelta_n\| = \max\{|t_i-t_{i-1}|: 1\le i\le n\}$ and assume $\|\itDelta_n\| \to
0$ as $n\to\infty$. Observe that
 \beqlb\label{2.13}
& &\E \<f, X^m_t\>
-
\E \bigg[\<Y_t, \mu^{M_t}\>
\exp\bigg\{\frac{1}{2}\int_0^t M_s(M_s-1)ds \bigg\}\bigg] \nonumber \\
&=&
\sum_{i=1}^n \bigg(\E \bigg[\<Y_{t-t_i},X_{t_i}^{M_{t-t_i}}\>
\exp\bigg\{\frac{1}{2}\int_0^{t-t_i} M_s(M_s-1)ds\bigg\}\bigg]  \\
& &\qquad
- \E \bigg[\<Y_{t-t_{i-1}},X_{t_{i-1}}^{M_{t-t_{i-1}}}\>
\exp\bigg\{\frac{1}{2}\int_0^{t-t_{i-1}} M_s(M_s-1)ds\bigg\}\bigg]\bigg). \nonumber
 \eeqlb
By the independence of $\{X_t: t\ge 0\}$ and $\{(M_t, Y_t): t\ge 0\}$ and the
martingale characterization of $\{(M_t, Y_t): t\ge 0\}$,
 \beqnn
& &\lim_{n\to \infty} \sum_{i=1}^n
\bigg(\E \bigg[\<Y_{t-t_i},X_{t_i}^{M_{t-t_i}}\>
\exp\bigg\{\frac{1}{2}\int_0^{t-t_i} M_s(M_s-1)ds\bigg\}\bigg]  \\
& &\qqquad
- \E \bigg[\<Y_{t-t_{i-1}},X_{t_i}^{M_{t-t_{i-1}}}\>
\exp\bigg\{\frac{1}{2}\int_0^{t-t_i} M_s(M_s-1)ds\bigg\}\bigg]\bigg)  \\
&=&
\lim_{n\to \infty} \sum_{i=1}^n
\E \bigg(\exp\bigg\{\frac{1}{2}\int_0^{t-t_i} M_s(M_s-1)ds\bigg\}
\E \bigg[F_{X_{t_i}}(M_{t-t_i},Y_{t-t_i})  \\
& &\qqquad
- F_{X_{t_i}}(M_{t-t_{i-1}},Y_{t-t_{i-1}})\bigg|X;\{(M_r,Y_r):0\le r\le
t-t_i\}\bigg]\bigg)   \\
&=&
- \lim_{n\to \infty} \sum_{i=1}^n
\E \bigg(\exp\bigg\{\frac{1}{2}\int_0^{t-t_i} M_s(M_s-1)ds\bigg\}  \\
& &\qqquad
\E \bigg[\int_{t-t_i}^{t-t_{i-1}} \L^*F_{X_{t_i}}(M_u,Y_u) du
\bigg| X;\{(M_r,Y_r):0\le r\le t-t_i\}\bigg]\bigg)   \\
&=&
- \lim_{n\to \infty} \sum_{i=1}^n
\E \bigg(\exp\bigg\{\frac{1}{2}\int_0^{t-t_i} M_s(M_s-1)ds\bigg\}
\int_{t-t_i}^{t-t_{i-1}} \L^*F_{X_{t_i}}(M_u,Y_u) du\bigg)   \\
&=&
- \lim_{n\to \infty} \int_0^t \sum_{i=1}^n
\E \bigg(\exp\bigg\{\frac{1}{2}\int_0^{t-t_i} M_s(M_s-1)ds\bigg\}
\L^*F_{X_{t_i}}(M_{t-u},Y_{t-u})\bigg) 1_{[t_{i-1},t_i]}(u) du   \\
&=&
- \int_0^t \E \bigg(\exp\bigg\{\frac{1}{2}\int_0^{t-u} M_s(M_s-1)ds\bigg\}
\L^*F_{X_u}(M_{t-u},Y_{t-u})\bigg) du,
 \eeqnn
where the last step holds by the right continuity of $\{X_t: t\ge 0\}$. Using again
the independence and the martingale problem for $\{X_t: t\ge 0\}$,
 \beqnn
& &\lim_{n\to \infty} \sum_{i=1}^n
\bigg(\E \bigg[\<Y_{t-t_{i-1}},X_{t_i}^{M_{t-t_{i-1}}}\>
\exp\bigg\{\frac{1}{2}\int_0^{t-t_i} M_s(M_s-1)ds\bigg\}\bigg]  \\
& &\qqquad
- \E \bigg[\<Y_{t-t_{i-1}},X_{t_{i-1}}^{M_{t-t_{i-1}}}\>
\exp\bigg\{\frac{1}{2}\int_0^{t-t_i} M_s(M_s-1)ds\bigg\}\bigg]\bigg)  \\
&=&
\lim_{n\to \infty} \sum_{i=1}^n
\E \bigg(\exp\bigg\{\frac{1}{2}\int_0^{t-t_i} M_s(M_s-1)ds\bigg\}  \\
& &\qqquad
\E \bigg[F_{M_{t-t_{i-1}},Y_{t-t_{i-1}}}(X_{t_i})
- F_{M_{t-t_{i-1}},Y_{t-t_{i-1}}}(X_{t_{i-1}})\bigg|M,Y\bigg]\bigg)  \\
&=&
\lim_{n\to \infty} \sum_{i=1}^n
\E \bigg(\exp\bigg\{\frac{1}{2}\int_0^{t-t_i} M_s(M_s-1)ds\bigg\}
\E \bigg[\int_{t_{i-1}}^{t_i}\L F_{M_{t-t_{i-1}},Y_{t-t_{i-1}}}(X_u)du
\bigg|M,Y\bigg]\bigg)  \\
&=&
\lim_{n\to \infty} \sum_{i=1}^n
\E \bigg(\exp\bigg\{\frac{1}{2}\int_0^{t-t_i} M_s(M_s-1)ds\bigg\}
\int_{t_{i-1}}^{t_i}\L F_{M_{t-t_{i-1}},Y_{t-t_{i-1}}}(X_u)du\bigg)  \\
&=&
\lim_{n\to \infty} \int_0^t \sum_{i=1}^n
\E \bigg(\exp\bigg\{\frac{1}{2}\int_0^{t-t_i} M_s(M_s-1)ds\bigg\}
\L F_{M_{t-t_{i-1}},Y_{t-t_{i-1}}}(X_u)\bigg) 1_{[t_{i-1},t_i]}(u) du  \\
&=&
\int_0^t
\E \bigg(\exp\bigg\{\frac{1}{2}\int_0^{t-u} M_s(M_s-1)ds\bigg\}
\L F_{M_{t-u},Y_{t-u}}(X_u)\bigg)du,
 \eeqnn
where we have also used the right continuity of $\{(M_t, Y_t): t\ge 0\}$ for the
last step. Finally, since $\|\itDelta_n\| \to 0$ as $n\to\infty$ and $M_t\le m$ for
all $t\ge0$, we have
 \beqnn
& &\lim_{n\to \infty} \sum_{i=1}^n
\bigg(\E \bigg[\<Y_{t-t_{i-1}},X_{t_{i-1}}^{M_{t-t_{i-1}}}\>
\exp\bigg\{\frac{1}{2}\int_0^{t-t_i} M_s(M_s-1)ds\bigg\}\bigg]  \\
& &\qqquad
- \E \bigg[\<Y_{t-t_{i-1}},X_{t_{i-1}}^{M_{t-t_{i-1}}}\>
\exp\bigg\{\frac{1}{2}\int_0^{t-t_{i-1}} M_s(M_s-1)ds\bigg\}\bigg]\bigg)  \\
&=&
\lim_{n\to \infty} \sum_{i=1}^n
\E \bigg(F_{X_{t_{i-1}}}(M_{t-t_{i-1}},Y_{t-t_{i-1}})
\exp\bigg\{\frac{1}{2}\int_0^{t-t_i} M_s(M_s-1)ds\bigg\}  \\
& &\qqquad
\bigg[1 - \exp\bigg\{\frac{1}{2}\int_{t-t_i}^{t-t_{i-1}}
M_u(M_u-1)du\bigg\}\bigg]\bigg)   \\
&=&
-\lim_{n\to \infty} \sum_{i=1}^n
\E \bigg(F_{X_{t_{i-1}}}(M_{t-t_{i-1}},Y_{t-t_{i-1}})
\exp\bigg\{\frac{1}{2}\int_0^{t-t_i} M_s(M_s-1)ds\bigg\}  \\
& &\qqquad
\bigg[\frac{1}{2}\int_{t-t_i}^{t-t_{i-1}} M_u(M_u-1)du\bigg]\bigg)   \\
&=&
-\lim_{n\to \infty} \frac{1}{2}\int_0^t \sum_{i=1}^n
\E \bigg(F_{X_{t_{i-1}}}(M_{t-t_{i-1}},Y_{t-t_{i-1}})
\exp\bigg\{\frac{1}{2}\int_0^{t-t_i} M_s(M_s-1)ds\bigg\}  \\
& &\qqquad
M_{t-u}(M_{t-u}-1)\bigg) 1_{[t_{i-1},t_i]}(u) du.
 \eeqnn
Since the semigroups $(P^m_t)_{t\ge0}$ are strongly Feller and strongly continuous,
$\{Y_t: t\ge 0\}$ is continuous in the uniform norm in each open interval between two
neighboring jumps of $\{M_t: t\ge 0\}$. Using this, the left continuity of $\{X_t:
t\ge 0\}$ and dominated convergence, we see that the above value is equal to
 \beqnn
- \frac{1}{2}\int_0^t
\E \bigg( F_{X_u}(M_{t-u},Y_{t-u})
\exp\bigg\{\frac{1}{2}\int_0^{t-u} M_s(M_s-1)ds\bigg\} M_{t-u}(M_{t-u}-1)\bigg)du.
 \eeqnn
Combining those together we see that the value of (\ref{2.13}) is in fact zero and
hence (\ref{2.11}) follows.
\qed

\btheorem\label{t2.2}
Let ${\cal D}(\L)$ be as in Theorem \ref{t2.1} and let $\{w_t: t\ge0\}$ denote the
coordinate process of $C([0,\infty), M(\IR))$. Suppose that for each $\mu \in
M(\IR)$ there is a probability measure $\Q_\mu$ on $C([0,\infty), M(\IR))$ such that
$\Q_\mu \{\<1,w_t\>^m\}$ is locally bounded in $t\ge0$ for every $m\ge1$ and such that
$\{w_t: t\ge0\}$ under $\Q_\mu$ is a solution of the $(\L, {\cal D}(\L),
\mu)$-martingale problem. Then the system $\{\Q_\mu: \mu \in M(\IR)\}$ defines a
diffusion process with transition semigroup $(Q_t)_{\ge0}$ given by
 \beqlb\label{2.14}
\int_{M(\IR)} \<f, \nu^m\>Q_t(\mu, d\nu)
=
\E^\sigma_{m,f} \bigg[\<Y_t, \mu^{M_t}\>
\exp\bigg\{\frac{1}{2}\int_0^t M_s(M_s-1)ds \bigg\}\bigg].
 \eeqlb
\etheorem

{\it Proof.} Let $Q_t(\mu, \cdot)$ denote the distribution of $w_t$ under
$\Q_\mu$. By Theorem \ref{t2.1} we have (\ref{2.14}). Let us assume first that
$\sigma(x) \equiv \sigma_0$ for a constant $\sigma_0$. In this case, $\{\<1,w_t\>:
t\ge0\}$ is the Feller diffusion with generator $(\sigma_0/2) x d^2/dx^2$, so that
 \beqnn
\int_{M(\IR)} e^{\lambda\<1,\nu\>} Q_t(\mu, d\nu)
=
\exp\bigg\{\frac{2\<1,\mu\>\lambda}{2-\sigma_0\lambda t}\bigg\},
\qquad t\ge0,\lambda\ge0.
 \eeqnn
Then for each $f\in B(\IR)^+$ the power series
 \beqlb\label{2.15}
\sum_{m=0}^\infty\frac{1}{m!}
\int_{M(\IR)} \<f,\nu\>^m Q_t(\mu, d\nu)\lambda^m
 \eeqlb
has a positive radius of convergence. By this and Billingsley (1968, p.342) it is not
hard to show that $Q_t(\mu, \cdot)$ is the unique probability measure on $M(\IR)$
satisfying (\ref{2.14}). Now the result follows from Ethier and Kurtz (1986, p.184).
For a non-constant $\sigma\in B(\IR)^+$, let $\sigma_0 = \|\sigma\|$ and observe that
 \beqnn
\int_{M(\IR)} \<f,\nu\>^m Q_t(\mu, d\nu)
\le
\E^{\sigma_0}_{m,f^{\otimes m}} \bigg[\<Y_t, \mu^{M_t}\>
\exp\bigg\{\frac{1}{2}\int_0^t M_s(M_s-1)ds \bigg\}\bigg]
 \eeqnn
by (\ref{2.14}) and the construction (\ref{2.6}) of $\{Y_t: t\ge0\}$, where $f^{\otimes
m}\in B(\IR^m)^+$ is defined by $f^{\otimes m}(x_1,\cdots,x_m) = f(x_1) \cdots
f(x_m)$. Then the power series (\ref{2.15}) also has a positive radius of convergence
and the result follows as in the case of a constant branching rate.
\qed


\section{Interacting-branching particle systems}
\setcounter{equation}{0}

In this section, we give a formulation of the interacting-branching particle system.
We first prove that equations (\ref{1.6}) have unique solutions. Recall that $c\in
C(\IR)$ is Lipschitz, $h\in C(\IR)$ is square-integrable and $\rho$ is twice
continuously differentiable with $\rho^{\prime}$ and $\rho^{\prime\prime}$ bounded.
The following result is an extension of Lemma 1.3 of Wang (1997) where it was assumed
that $c(x)\equiv$ const.

\blemma\label{l3.1}
For any initial conditions $x_i(0) = x_i$, equations (\ref{1.6}) have unique
solutions $\{x_i(t): t\ge 0\}$ and $\{(x_1(t), \cdots, x_m(t)): t\ge0\}$ is an
$m$-dimensional diffusion process with generator $G^m$ defined by (\ref{1.7}).
\elemma

{\it Proof.} Fix $T>0$ and $i\ge1$ and define $\{x_i^k(t): t\ge 0\}$ inductively by
$x_i^0(t) \equiv x_i(0)$ and
 \beqnn
x_i^{k+1}(t)
=
x_i(0) + \int_0^tc(x_i^k(s)) dB_i(s)
+ \int_0^t\int_{\IR} h(y-x_i^k(s))W(dy, ds),
\qquad t\ge0.
 \eeqnn
Let $l(c)\ge 0$ be any Lipschitz constant for $c(\cdot)$. By a martingale inequality we
have
 \beqnn
\E\bigg\{\sup_{0\le t\le T}|x_i^{k+1}(t) - x_i^k(t)|^2\bigg\}
&\le&
8\int_0^T \E\{|c(x_i^k(t)) - c(x_i^{k-1}(t))|^2\} dt  \\
& &
+ 8\int_0^T \E\bigg\{\int_{\IR}|h(y-x_i^k(t))
- h(y-x_i^{k-1}(t))|^2 dy\bigg\} dt  \\
&\le&
8l(c)^2\int_0^T\E\{|x_i^k(t) - x_i^{k-1}(t)|^2\} dt  \\
& & +
16 \int_0^T \E\{|\rho(0) - \rho(x_i^k(t) - x_i^{k-1}(t))|\} dt  \\
&\le&
8 (l(c)^2 + \|\rho^{\prime\prime}\|)
\int_0^T \E\{|x_i^k(t) - x_i^{k-1}(t)|^2\} dt.
 \eeqnn
Using the above inequality inductively we get
 \beqnn
\E\bigg\{\sup_{0\le t\le T}|x_i^{k+1}(t) - x_i^k(t)|^2\bigg\}
\le
(\|c\|^2 +\rho(0)) (l(c)^2 + \|\rho^{\prime\prime}\|)^k (8T)^k/k!,
 \eeqnn
and hence
 \beqnn
\P\bigg\{\sup_{0\le t\le T}|x_i^{k+1}(t) - x_i^k(t)| > 2^{-k}\bigg\}
\le
\mbox{const}\cdot (l(c)^2 + \|\rho^{\prime\prime}\|)^k(8T)^k/k!.
 \eeqnn
By Borel-Cantelli's lemma, $\{x_i^k(t): 0\le t\le T\}$ converges in the uniform
norm with probability one. Since $T>0$ was arbitrary, $x_i(t) = \lim_{k\to\infty}
x_i^k(t)$ defines a continuous martingale $\{x_i(t): t\ge 0\}$ which is clearly the
unique solution of (\ref{1.6}). It is easy to see that $d\<x_i\>(t) = a(x_i(t))dt$ and
$d\<x_i,x_j\>(t) = \rho(x_i(t) - x_j(t))dt$ for $i\neq j$. Then $\{(x_1(t), \cdots,
x_m(t)): t\ge0\}$ is a diffusion process with generator $G^m$ defined by (\ref{1.7}).
\qed

Because of the exchangebility, the $G^m$-diffusion can be regarded as a measure-valued
Markov process. Let $N(\IR)$ denote the space of integer-valued measures on $\IR$.
For $\theta >0$, let $M_\theta(\IR) = \{\theta^{-1} \sigma: \sigma\in N(\IR)\}$. Let
$\zeta$ be the mapping from $\cup_{m=1} ^\infty\IR^m$ to $M_\theta(\IR)$ defined by
 \beqlb\label{3.1}
\zeta(x_1,\cdots,x_m) = \frac{1}{\theta} \sum_{i=1}^m\delta_{x_i},
\qquad m\ge1.
 \eeqlb

\blemma\label{l3.2}
For any integers $m,n\ge1$ and any $f\in C^2 (\IR^n)$, we have
 \beqlb\label{3.2}
G^m F_{n,f}(\zeta(x_1,\cdots,x_m))
&=&
\frac{1}{2\theta^n}\sum_{\alpha=1}^n
\sum_{l_1,\cdots,l_n=1}^m a(x_{l_\alpha})
f^{\prime\prime}_{\alpha\alpha} (x_{l_1},\cdots,x_{l_n})   \nonumber  \\
&+&
\frac{1}{2\theta^n}\sum_{\alpha,\beta=1, \alpha \neq \beta}^n
\sum_{l_1,\cdots,l_n=1, l_\alpha = l_\beta}^m c(x_{l_\alpha}) c(x_{l_\beta})
f^{\prime\prime}_{\alpha\beta} (x_{l_1},\cdots,x_{l_n})  \nonumber \\
&+&
\frac{1}{2\theta^n}\sum_{\alpha,\beta=1, \alpha \neq \beta}^n
\sum_{l_1,\cdots,l_n=1}^m
\rho(x_{l_\alpha}-x_{l_\beta})f^{\prime\prime}_{\alpha\beta}
(x_{l_1},\cdots,x_{l_n}).
 \eeqlb
\elemma

{\it Proof.} By (\ref{3.1}), we have
 \beqlb\label{3.3}
F_{n,f}(\zeta(x_1,\cdots,x_m))
=
\frac{1}{\theta^n}\sum_{l_1,\cdots,l_n=1}^mf(x_{l_1},\cdots,x_{l_n}).
 \eeqlb
Observe that, for $1\le i\le m$,
 \beqnn
\frac{d^2}{dx_i^2}F_{n,f}(\zeta(x_1,\cdots,x_m))
=
\frac{1}{\theta^n}\sum_{\alpha,\beta=1}^n\sum_{\{\cdots\}}
f^{\prime\prime}_{\alpha\beta} (x_{l_1},\cdots,x_{l_n}),
 \eeqnn
where $\{\cdots\} =\{$\,for all $1\le l_1,\cdots,l_n \le m$ with $l_\alpha = l_\beta
=i\}$. Then it is not hard to see that
 \beqlb\label{3.4}
& & \sum_{i=1}^m c(x_i)^2 \frac{d^2}{dx_i^2}F_{n,f}
(\zeta(x_1,\cdots,x_m)) \nonumber \\
&=&
\frac{1}{\theta^n}\sum_{\alpha,\beta=1}^n
\sum_{l_1,\cdots,l_n=1, l_\alpha = l_\beta}^m c(x_{l_\alpha}) c(x_{l_\beta})
f^{\prime\prime}_{\alpha\beta} (x_{l_1},\cdots,x_{l_n})  \nonumber \\
&=&
\frac{1}{\theta^n}\sum_{\alpha=1}^n
\sum_{l_1,\cdots,l_n=1}^m c(x_{l_\alpha})^2
f^{\prime\prime}_{\alpha\alpha} (x_{l_1},\cdots,x_{l_n})  \nonumber \\
& &\quad + \frac{1}{\theta^n}\sum_{\alpha,\beta=1, \alpha \neq \beta}^n
\sum_{l_1,\cdots,l_n=1, l_\alpha = l_\beta}^m c(x_{l_\alpha}) c(x_{l_\beta})
f^{\prime\prime}_{\alpha\beta} (x_{l_1},\cdots,x_{l_n}).
 \eeqlb
On the other hand, for $1\le i \neq j\le m$,
 \beqnn
\bigg(\frac{d^2}{dx_idx_j} + \frac{d^2}{dx_idx_j}\bigg)
F_{n,f}(\zeta(x_1,\cdots,x_m))
=
\frac{1}{\theta^n}\sum_{\alpha,\beta=1,\alpha \neq \beta}^n\sum_{\{\cdots\}}
f^{\prime\prime}_{\alpha\beta} (x_{l_1},\cdots,x_{l_n}),
 \eeqnn
where $\{\cdots\} =\{$\,for all $1\le l_1,\cdots,l_n \le m$ with $l_\alpha = i$ and
$l_\beta =j\}$. It follows that
 \beqnn
& & \sum_{i,j=1,i\neq j}^m\rho(x_i-x_j)\frac{d^2}{dx_idx_j}
F_{n,f}(\zeta(x_1,\cdots,x_m))  \\
&=&
\frac{1}{\theta^n}\sum_{\alpha,\beta=1, \alpha \neq \beta}^n
\sum_{l_1,\cdots,l_n=1, l_\alpha \neq l_\beta}^m
\rho(x_{l_\alpha}-x_{l_\beta}) f^{\prime\prime}_{\alpha\beta}
(x_{l_1},\cdots,x_{l_n}).
 \eeqnn
Using this and (\ref{3.4}) with $c(x_i)^2$ replaced by $\rho(0)$,
 \beqlb\label{3.5}
& & \sum_{i,j=1}^m\rho(x_i-x_j)\frac{d^2}{dx_idx_j}
F_{n,f}(\zeta(x_1,\cdots,x_m))  \nonumber \\
&=&
\frac{1}{\theta^n} \sum_{\alpha=1}^n
\sum_{l_1,\cdots,l_n=1}^m \rho(0)
f^{\prime\prime}_{\alpha\alpha} (x_{l_1},\cdots,x_{l_n}) \nonumber  \\
& &+\,\frac{1}{\theta^n}\sum_{\alpha,\beta=1, \alpha \neq \beta}^n
\sum_{l_1,\cdots,l_n=1}^m
\rho(x_{l_\alpha}-x_{l_\beta})f^{\prime\prime}_{\alpha\beta}
(x_{l_1},\cdots,x_{l_n}).
 \eeqlb
Then we have the desired result from (\ref{3.4}) and (\ref{3.5}). \qed

Suppose that $X(t) = (x_1(t),\cdots,x_m(t))$ is a Markov process in $\IR^m$ generated
by $G^m$. Based on (\ref{1.2}) and Lemma \ref{l3.2}, it is easy to show that $\zeta
(X(t))$ is a Markov process in $M_\theta(\IR)$ with generator $\A_\theta$ given by
 \beqlb\label{3.6}
\A_\theta F(\mu)
&=&
\frac{1}{2}\int_{\IR} a(x)\frac{d^2}{dx^2}\frac{\delta
F(\mu)}{\delta\mu(x)}\mu(dx)
+ \frac{1}{2\theta}\int_{\IR^2} c(x)c(y) \frac{d^2}{dxdy}
\frac{\delta^2 F(\mu)}{\delta\mu(x)\delta\mu(y)}\delta_x(dy)\mu(dx) \nonumber  \\
& & +\,\frac{1}{2}\int_{\IR^2}\rho(x-y)\frac{d^2}{dxdy}
\frac{\delta^2 F(\mu)}{\delta\mu(x)\delta\mu(y)}\mu(dx)\mu(dy).
 \eeqlb
In particular, if
 \beqlb\label{3.7}
F(\mu)
= f(\<\phi_1,\mu\>, \cdots, \<\phi_n,\mu\>),
\qquad \mu\in M_\theta(\IR),
 \eeqlb
for $f\in C^2(\IR^n)$ and $\{\phi_i\}\subset C^2 (\IR)$, then
 \beqlb\label{3.8}
\A_\theta F(\mu)
&=&
\frac{1}{2}\sum_{i=1}^n f_i^\prime(\<\phi_1,\mu\>, \cdots, \<\phi_n,\mu\>)
\<a\phi_i^{\prime\prime},\mu\>       \nonumber \\
&+&
\frac{1}{2\theta}\sum_{i,j=1}^n f_{ij}^{\prime\prime}(\<\phi_1,\mu\>, \cdots,
\<\phi_n,\mu\>) \<c^2\phi_i^\prime\phi_j^\prime,\mu\>  \\
&+&
\frac{1}{2}\sum_{i,j=1}^n f_{ij}^{\prime\prime}(\<\phi_1,\mu\>, \cdots,
\<\phi_n,\mu\>)\int_{\IR^2}\rho(x-y) \phi_i^\prime(x) \phi_j^\prime(y)
\mu(dx)\mu(dy). \nonumber
 \eeqlb

Now we introduce a branching mechanism to the interacting particle system. Suppose
that for each $x\in \IR$ we have a discrete probability distribution $p(x) = \{p_i(x):
i=0,1,\cdots\}$ such that each $p_i(\cdot)$ is a Borel measurable function on $\IR$.
This serves as the distribution of the offspring number produced by a particle that
dies at site $x\in \IR$. We assume that
 \beqlb\label{3.9}
p_1(x) = 0,
\quad
\sum_{i=1}^\infty ip_i(x) =1,
 \eeqlb
and
 \beqlb\label{3.10}
\sigma_p(x) := \sum_{i=1}^\infty i^2p_i(x) -1
 \eeqlb
is bounded in $x\in \IR$. Let $\itGamma_\theta (\mu,d\nu)$ be the probability kernel
on $M_\theta(\IR)$ defined by
 \beqlb\label{3.11}
\int_{M_\theta(\IR)} F(\nu) \itGamma_\theta (\mu,d\nu)
=
\frac{1}{\theta\mu(1)}\sum_{i=1}^{\theta\mu(1)}\sum_{j=0}^\infty
p_j(x_i) F\bigg(\mu + (j-1)\theta^{-1} \delta_{x_i}\bigg),
 \eeqlb
where $\mu\in M_\theta(\IR)$ is given by
 \beqnn
\mu = \frac{1}{\theta} \sum_{i=1}^{\theta\mu(1)} \delta_{x_i}.
 \eeqnn
For a constant $\gamma>0$, we define the bounded operator ${\cal B}_\theta$ on
$B(M_\theta(\IR))$ by
 \beqlb\label{3.12}
\B_\theta F(\mu)
=
\gamma\theta^2 [\theta\land \mu(1)] \int_{M_\theta(\IR)}
[F(\nu) - F(\mu)] \itGamma_\theta (\mu,d\nu).
 \eeqlb
In view of (\ref{1.6}), $\A_\theta$ generates a Feller Markov process on
$M_\theta(\IR)$, then so does $\L_\theta := \A_\theta + \B_\theta$ by Ethier-Kurtz
(1986, p.37). We shall call the process generated by $\L_\theta$ an {\it
interacting-branching particle system} with parameters $(a,\rho,\gamma,p)$ and unit
mass $1/\theta$. Heuristically, each particle in the system has mass $1/\theta$,
$a(\cdot)$ represents the migration speed of the particles and $\rho(\cdot)$
describes the interaction between them. The branching times of the system are
determined by the killing density $\gamma\theta^2 [\theta\land \mu(1)]$, where the
truncation ``$\theta\land \mu(1)$'' is introduced to make the branching not too fast
even when the total mass is large. At each branching time, with equal probability, one
particle in the system is randomly chosen, which is killed at its site $x\in \IR$ and
the offspring are produced at $x\in \IR$ according to the distribution $\{p_i(x):
i=0,1,\cdots\}$. If $F$ is given by (\ref{3.7}), then $\B_\theta F(\mu)$ is equal to
 \beqlb\label{3.13}
\frac{\gamma [\theta\land \mu(1)]}{2\mu(1)}\sum_{\alpha,\beta=1}^n
\sum_{j=1}^\infty (j-1)^2
\<p_jf_{\alpha\beta}^{\prime\prime}(\<\phi_1,\mu\> +\xi_j\phi_1, \cdots,
\<\phi_n,\mu\> +\xi_j\phi_n) \phi_\alpha\phi_\beta,\mu\>
 \eeqlb
for some constant $0< \xi_j <(j-1)/\theta$. This follows from (\ref{3.11}) and
(\ref{3.12}) by Taylor's expansion.


\section{Continuous branching density}
\setcounter{equation}{0}
In this section, we shall construct a solution of the martingale problem of the SDSM
with continuous branching density by using particle system approximation. Assume
that $\sigma \in C(\IR)$ can be extended continuously to $\hat \IR$. Let $\A$ and
$\B$ be given by (\ref{1.5}) and (\ref{1.8}), respectively. Observe that, if
 \beqlb\label{4.1}
F(\mu) = f(\<\phi_1,\mu\>, \cdots, \<\phi_n,\mu\>),
\qquad \mu\in M(\IR),
 \eeqlb
for $f\in C^2(\IR^n)$ and $\{\phi_i\}\subset C^2(\IR)$, then
 \beqlb\label{4.2}
\A F(\mu)
&=&
\frac{1}{2}\sum_{i=1}^n f_i^\prime(\<\phi_1,\mu\>, \cdots,
\<\phi_n,\mu\>) \<a\phi_i^{\prime\prime},\mu\>  \\
&+&
\frac{1}{2}\sum_{i,j=1}^n
f_{ij}^{\prime\prime}(\<\phi_1,\mu\>, \cdots, \<\phi_n,\mu\>) \int_{\IR^2}
\rho(x-y) \phi_i^\prime(x) \phi_j^\prime(y) \mu(dx)\mu(dy),  \nonumber
 \eeqlb
and
 \beqlb\label{4.3}
\B F(\mu)
&=&
\frac{1}{2} \sum_{i,j=1}^n
f_{ij}^{\prime\prime}(\<\phi_1,\mu\>, \cdots, \<\phi_n,\mu\>)
\<\sigma \phi_i\phi_j,\mu\>.
 \eeqlb

Let $\{\theta_k\}$ be any sequence such that $\theta_k\to \infty$ as $k\to \infty$.
Suppose that $\{X_t^{(k)}: t\ge0\}$ is a sequence of c\'adl\'ag interacting-branching
particle systems with parameters $(a, \rho, \gamma_k, p^{(k)})$, unit mass
$1/\theta_k$ and initial states $X_0^{(k)} = \mu_k \in M_{\theta_k}(\IR)$. In an
obvious way, we may also regard $\{X_t^{(k)}: t\ge0\}$ as a process with state space
$M(\hat\IR)$. Let $\sigma_k$ be defined by (\ref{3.10}) with $p_i$ replaced by
$p_i^{(k)}$.

\blemma\label{l4.1}
Suppose that the sequences $\{\gamma_k \sigma_k\}$ and $\{\<1, \mu_k\>\}$ are
bounded. Then $\{X_t^{(k)}: t\ge0\}$ form a tight sequence in $D([0,\infty),
M(\hat\IR))$. \elemma

{\it Proof.} By the assumption (\ref{3.9}), it is easy to show that $\{\<1, X_t^{(k)}\>:
t\ge0\}$ is a martingale. Then we have
 \beqnn
\P\bigg\{\sup_{t\ge0}\<1, X_t^{(k)}\> > \eta\bigg\}
\le
\frac{\<1, \mu_k\>}{\eta}
 \eeqnn
for any $\eta>0$. That is, $\{X_t^{(k)}: t\ge0\}$ satisfies the compact containment
condition of Ethier and Kurtz (1986, p.142). Let $\L_k$ denote the generator of
$\{X_t^{(k)}: t\ge0\}$ and let $F$ be given by (\ref{4.1}) with $f\in C^2_0 (\IR^n)$
and with each $\phi_i\in C^2_\partial (\IR)$ bounded away from zero. Then
 \beqnn
F(X_t^{(k)}) - F(X_0^{(k)}) - \int_0^t \L_k F(X_s^{(k)})ds,
\qquad t\ge0,
 \eeqnn
is a martingale and the desired tightness follows from the result of Ethier and Kurtz
(1986, p.145).
\qed

In the sequel of this section, we assume $\{\phi_i\}\subset C^2_\partial (\IR)$.
In this case, (\ref{4.1}), (\ref{4.2}) and (\ref{4.3}) can be extended to
continuous functions on $M(\hat\IR)$. Let
$\hat\A F(\mu)$ and $\hat\B F(\mu)$ be defined respectively by the right hand side of
(\ref{4.2}) and (\ref{4.3}) and let $\hat\L F(\mu) = \hat\A F(\mu) + \hat\B
F(\mu)$, all defined as continuous functions on $M(\hat\IR)$.

\blemma\label{l4.2}
Let ${\cal D}(\hat\L)$ be the totality of all functions of the form (\ref{4.1}) with
$f\in C^2_0 (\IR^n)$ and with each $\phi_i\in C^2_\partial (\IR)$ bounded away
from zero. Suppose further that $\gamma_k \sigma_k \to \sigma$ uniformly and $\mu_k
\to \mu \in M(\hat \IR)$ as $k\to \infty$. Then any limit point $\Q_\mu$ of the
distributions of $\{X_t^{(k)}: t\ge0\}$ is supported by $C([0,\infty), M(\hat\IR))$
under which
 \beqlb\label{4.4}
F(w_t) - F(w_0) - \int_0^t\hat\L F(w_s)ds,
\qquad t\ge 0,
 \eeqlb
is a martingale for each $F\in {\cal D}(\hat\L)$, where $\{w_t: t\ge0\}$ denotes the
coordinate process of $C([0,\infty), M(\hat\IR))$. \elemma

{\it Proof.} We use the notation introduced in the proof of Lemma
\ref{l4.1}. By passing to a subsequence if it is necessary, we may
assume that the distribution of $\{X_t^{(k)}: t\ge0\}$ on
$D([0,\infty), M(\hat\IR))$ converges to $\Q_\mu$. Using
Skorokhod's representation, we may assume that the processes
$\{X_t^{(k)}: t\ge0\}$ are defined on the same probability space
and the sequence converges almost surely to a c\'adl\'ag process
$\{X_t: t\ge0\}$ with distribution $\Q_\mu$ on $D([0,\infty),
M(\hat\IR))$; see e.g. Ethier and Kurtz (1986, p.102). Let $K(X) =
\{t\ge 0: \P\{X_t = X_{t^-}\} =1\}$. By Ethier and Kurtz (1986,
p.118), for each $t\in K(X)$ we have a.s. $\lim_{k\to \infty}
X_t^{(k)} =  X_t$. Recall that $f$ and $f^{\prime\prime}_{ij}$ are
rapidly decreasing and each $\phi_i$ is bounded away from zero.
Since $\gamma_ka_k \to \sigma$ uniformly, for $t\in K(X)$ we have
a.s. $\lim_{k\to \infty} \L_k F(X_t^{(k)}) = \hat\L F(X_t)$
boundedly by (\ref{3.8}), (\ref{3.13}) and the definition of
$\hat\L$. Suppose that $\{H_i\}_{i=1}^n \subset C(M(\hat\IR))$ and
$\{t_i\}_{i=1}^{n+1} \subset K(X)$ with $0\le t_1 <\cdots <t_n
<t_{n+1}$. By Ethier and Kurtz (1986, p.131), the complement of
$K(X)$ is at most countable. Then
 \beqnn
& &\E\bigg\{\bigg[F(X_{t_{n+1}}) - F(X_{t_n})
- \int_{t_n}^{t_{n+1}}\hat\L F(X_s)ds\bigg]
\prod_{i=1}^n H_i(X_{t_i})\bigg\}  \\
&=&
\E\bigg\{F(X_{t_{n+1}})\prod_{i=1}^n H_i(X_{t_i})\bigg\}
- \E\bigg\{F(X_{t_n})\prod_{i=1}^n H_i(X_{t_i})\bigg\} \\
& &\qqquad
- \int_{t_n}^{t_{n+1}}\E\bigg\{\hat\L F(X_s)\prod_{i=1}^n H_i(X_{t_i})\bigg\}ds  \\
&=&
\lim_{k\to \infty}\E\bigg\{F(X^{(k)}_{t_{n+1}})
\prod_{i=1}^n H_i(X^{(k)}_{t_i})\bigg\}
- \lim_{k\to \infty}\E\bigg\{F(X^{(k)}_{t_n})
\prod_{i=1}^n H_i(X^{(k)}_{t_i})\bigg\} \\
& &\qqquad - \lim_{k\to \infty}\int_{t_n}^{t_{n+1}}\E\bigg\{\L_k F(X^{(k)}_s)
\prod_{i=1}^n H_i(X^{(k)}_{t_i})\bigg\}ds  \\
&=&
\lim_{k\to \infty}\E\bigg\{\bigg[F(X^{(k)}_{t_{n+1}}) - F(X^{(k)}_{t_n}) -
\int_{t_n}^{t_{n+1}}\L_k F(X^{(k)}_s)ds\bigg]
\prod_{i=1}^n H_i(X^{(k)}_{t_i})\bigg\} \\
&=& 0.
 \eeqnn
By the right continuity of $\{X_t: t\ge0\}$, the equality
 \beqnn
\E\bigg\{\bigg[F(X_{t_{n+1}}) - F(X_{t_n}) - \int_{t_n}^{t_{n+1}}
\hat\L F(X_s)ds\bigg] \prod_{i=1}^n H_i(X_{t_i})\bigg\}
= 0
 \eeqnn
holds without the restriction $\{t_i\}_{i=1}^{n+1} \subset K(X)$. That is,
 \beqnn
F(X_t) - F(X_0) - \int_0^t\hat\L F(X_s)ds,
\qquad t\ge 0,
 \eeqnn
is a martingale. As in Wang (1998, pp.783-784) one can show that $\{X_t: t\ge0\}$ is
in fact a.s. continuous.
\qed

\blemma\label{l4.3}
Let ${\cal D}(\hat\L)$ be as in Lemma \ref{l4.2}. Then for each $\mu \in M(\hat\IR)$,
there is a probability measure $\Q_\mu$ on $C([0,\infty), M(\hat\IR))$ under which
(\ref{4.4}) is a martingale for each $F\in {\cal D}(\hat\L)$.
\elemma

{\it Proof.} It is easy to find $\mu_k \in M_{\theta_k} (\IR)$ such that $\mu_k \to
\mu$ as $k\to \infty$. Then, by Lemma \ref{l4.2}, it suffices to construct a sequence
$(\gamma_k,p^{(k)})$ such that $\gamma_k \sigma_k \to \sigma$ as $k\to \infty$. This
is elementary. One choice is described as follows. Let $\gamma_k = 1/\sqrt{k}$ and
$\sigma_k = \sqrt{k}(\sigma + 1/\sqrt{k})$. Then the system of equations
 \beqnn
\left\{\begin{array}{ll}
p_0^{(k)} + p_2^{(k)}  + p_k^{(k)}    &= 1,  \\
           2p_2^{(k)} + kp_k^{(k)}    &= 1,  \\
           4p_2^{(k)} + k^2p_k^{(k)}  &= \sigma_k +1,
\end{array}\right.
 \eeqnn
has the unique solution
 \beqnn
p_0^{(k)} = \frac{\sigma_k+k-1}{2k},
\quad
p_2^{(k)} = \frac{k-1-\sigma_k}{2(k-2)},
\quad
p_k^{(k)} = \frac{\sigma_k-1}{k(k-2)},
 \eeqnn
where each $p_i^{(k)}$ is nonnegative for sufficiently large $k\ge3$. \qed

\blemma\label{l4.4}
Let $\Q_\mu$ be given by Lemma \ref{l4.3}. Then for $n\ge1$, $t\ge0$ and $\mu\in
M(\IR)$ we have
 \beqnn
\Q_\mu\{\<1,w_t\>^n\}
\le
\<1,\mu\>^n + \frac{1}{2} n(n-1)\|\sigma\|
\int_0^t\Q_\mu\{\<1,w_s\>^{n-1}\}ds.
 \eeqnn
Consequently, $\Q_\mu\{\<1,w_t\>^n\}$ is a locally bounded function of $t\ge0$. Let
${\cal D}(\hat\L)$ be the union of all functions of the form (\ref{4.1}) with $f\in
C^2_0 (\IR^n)$ and $\{\phi_i\}\subset C^2_\partial (\IR)$ and all functions of the form
$F_{m,f}(\mu) = \<f,\mu^m\>$ with $f\in C^2_\partial (\IR^m)$. Then (\ref{4.4})
under $\Q_\mu$ is a martingale for each $F\in {\cal D}(\hat\L)$.
\elemma

{\it Proof.} For any $k\ge1$, take $f_k\in C^2_0(\IR))$ such that $f_k(z) = z^n$ for
$0\le z\le k$ and $f_k^{\prime\prime}(z) \le n(n-1)z^{n-2}$ for all $z\ge 0$. Let
$F_k(\mu) = f_k(\<1,\mu\>)$. Then $\A F_n(\mu) =0$ and
 \beqnn
\B F_k(\mu)
\le
\frac{1}{2} n(n-1) \|\sigma\| \<1,\mu\>^{n-1}.
 \eeqnn
Since
 \beqnn
F_k(X_t) - F_k(X_0) - \int_0^t\L F_k(\<1,X_s\>)ds,
\qquad t\ge0,
 \eeqnn
is a martingale, we get
 \beqnn
\Q_\mu f_k(\<1,X_t\>^n)
&\le&
f_k(\<1,\mu\>)
+ \frac{1}{2} n(n-1) \|\sigma\|
\int_0^t \Q_\mu(\<1,X_s\>^{n-1})ds  \\
&\le&
\<1,\mu\>^n
+ \frac{1}{2} n(n-1) \|\sigma\|
\int_0^t \Q_\mu(\<1,X_s\>^{n-1})ds.
 \eeqnn
Then the desired estimate follows by Fatou's Lemma. The last assertion is an
immediate consequence of Lemma \ref{l4.3}.
\qed

\blemma\label{l4.5}
Let $\Q_\mu$ be given by Lemma \ref{l4.3}. Then for $\mu \in M(\IR)$ and $\phi \in
C^2_\partial(\IR)$,
 \beqlb\label{4.5}
M_t(\phi)
:=
\<\phi, w_t\> - \<\phi, \mu\>
- \frac{1}{2} \int_0^t \<a\phi^{\prime\prime}, w_s\> ds,
\qquad t\ge0,
 \eeqlb
is a $\Q_\mu$-martingale with quadratic variation process
 \beqlb\label{4.6}
\<M(\phi)\>_t
=
\int_0^t\<\sigma\phi^2, w_s\>ds +
\int_0^t ds\int_{\hat\IR} \<h(z - \cdot) \phi^\prime, w_s\>^2 dz.
 \eeqlb
\elemma

{\it Proof.} It is easy to check that, if $F_n(\mu) = \<\phi,\mu\>^n$, then
 \beqnn
\hat\L F_n(\mu)
&=&
\frac{n}{2} \<\phi,\mu\>^{n-1} \<a\phi^{\prime\prime},\mu\>
+ \frac{n(n-1)}{2} \<\phi,\mu\>^{n-2}
\int_{\hat\IR} \<h(z-\cdot)\phi^{\prime}, \mu\>^2 dz         \\
& &\qqquad\qquad + \frac{n(n-1)}{2} \<\phi,\mu\>^{n-2}\<\sigma\phi^2,\mu\>.
 \eeqnn
It follows that both (\ref{4.5}) and
 \beqlb\label{4.7}
M_t^2(\phi)
&:=&
\<\phi, w_t\>^2 - \<\phi, \mu\>^2 - \int_0^t \<\phi, w_s\>
\<a\phi^{\prime\prime}, w_s\>ds    \nonumber \\
& &\qqquad
- \int_0^tds\int_{\hat\IR} \<h(z-\cdot)\phi^\prime, w_s\>^2 dz
- \int_0^t \<\sigma\phi^2, w_s\>ds
 \eeqlb
are martingales. By (\ref{4.5}) and It\^o's formula we have
 \beqlb\label{4.8}
\<\phi, w_t\>^2
=
\<\phi, \mu\>^2
+ \int_0^t \<\phi, w_s\> \<a\phi^{\prime\prime}, w_s\>ds
+ 2\int_0^t \<\phi, w_s\> dM_s(\phi) + \<M(\phi)\>_t.
 \eeqlb
Comparing (\ref{4.7}) and (\ref{4.8}) we get the conclusion. \qed

Observe that the martingales $\{M_t(\phi): t\ge 0\}$ defined by (\ref{4.5}) form a system
which is linear in $\phi \in C^2_\partial (\IR)$. Because of the presence of the
derivative $\phi^\prime$ in the variation process (\ref{4.6}), it seems hard to extend the
definition of $\{M_t(\phi): t\ge 0\}$ to a general function $\phi \in B(\hat\IR)$.
However, following the method of Walsh (1986), one can still define the stochastic
integral
 \beqnn
\int_0^t\int_{\hat\IR}\phi(s,x)M(ds,dx),
\qquad t\ge 0,
 \eeqnn
if both $\phi(s,x)$ and $\phi^\prime_x(s,x)$ can be extended continuously to
$[0,\infty) \times \hat\IR$. With those in hand, we have the following

\blemma\label{l4.6}
Let $\Q_\mu$ be given by Lemma \ref{l4.3}. Then for any $t\ge0$ and $\phi \in
C^2_\partial (\IR)$ we have a.s.
 \beqnn
\<\phi, w_t\>
=
\<\hat P_t\phi, \mu\> + \int_0^t\int_{\hat\IR}\hat P_{t-s}\phi(x)
M(ds,dx).
 \eeqnn
\elemma

{\it Proof.} For any partition $\itDelta_n := \{0= t_0 <t_1 <\cdots <t_n =t\}$ of
$[0,t]$, we have
 \beqnn
\<\phi,w_t\> - \<\hat P_t\phi,\mu\>
&=&
\sum_{i=1}^n \<\hat P_{t-t_i}\phi - \hat P_{t-t_{i-1}}\phi, w_{t_i}\> \\
& &\quad
+ \sum_{i=1}^n [\<\hat P_{t-t_{i-1}}\phi,w_{t_i}\>
- \<\hat P_{t-t_{i-1}}\phi,w_{t_{i-1}}\>].
 \eeqnn
Let $\|\itDelta_n\| = \max\{|t_i-t_{i-1}|: 1\le i\le n\}$ and assume $\|\itDelta_n\|
\to 0$ as $n\to\infty$. Then
 \beqnn
\lim_{n\to \infty} \sum_{i=1}^n \<\hat P_{t-t_i}\phi
-\hat P_{t-t_{i-1}}\phi,w_{t_i}\>
&=&
- \lim_{n\to \infty} \sum_{i=1}^n \int_{t_{i-1}}^{t_i}
\<\hat P_{t-s}\hat G\phi, w_{t_i}\>ds \\
&=&
-\int_0^t\<\hat P_{t-s} \hat G \phi, w_s\>ds.
 \eeqnn
Using Lemma \ref{l4.5} we have
 \beqnn
& &\lim_{n\to \infty} \sum_{i=1}^n [\<\hat P_{t-t_{i-1}}\phi,w_{t_i}\>
- \<\hat P_{t-t_{i-1}}\phi,w_{t_{i-1}}\>]    \\
&=&
\lim_{n\to \infty} \sum_{i=1}^n \int_{t_{i-1}}^{t_i}\int_{\hat\IR}
\hat P_{t-t_{i-1}}\phi(x) M(ds,dx)
+ \lim_{n\to \infty} \frac{1}{2}\sum_{i=1}^n
\int_{t_{i-1}}^{t_i}\<a(\hat P_{t-t_{i-1}}\phi)^{\prime\prime}, w_s\> ds \\
&=&
\int_0^t\int_{\hat\IR}\hat P_{t-s}\phi(x) M(ds,dx)
+ \frac{1}{2}\int_0^t \<a(\hat P_{t-s}\phi)^{\prime\prime}, w_s\> ds.
 \eeqnn
Combining those we get the desired conclusion.
\qed

\btheorem\label{t4.1}
Let ${\cal D}(\L)$ be the union of all functions of the form (\ref{4.1}) with $f\in
C^2 (\IR^n)$ and $\{\phi_i\}\subset C^2 (\IR)$ and all functions of the form $F_{m,f} (\mu)
= \<f,\mu^m\>$ with $f\in C^2 (\IR^m)$. Let $\{w_t: t\ge0\}$ denote the coordinate
process of $C([0,\infty), M(\IR))$. Then for each $\mu \in M(\IR)$ there is a
probability measure $\Q_\mu$ on $C([0,\infty), M(\IR))$ such that $\Q_\mu
\{\<1,w_t\>^m\}$ is locally bounded in $t\ge0$ for every $m\ge1$ and such that
$\{w_t: t\ge0\}$ under $\Q_\mu$ is a solution of the $(\L, {\cal D}(\L),
\mu)$-martingale problem. \etheorem

{\it Proof.} Let $\Q_\mu$ be the probability measure on $C([0,\infty), M(\hat\IR))$
provided by Lemma \ref{l4.3}. The desired result will follow once it is proved that
 \beqlb\label{4.9}
\Q_\mu\{w_t(\{\partial\})=0 \mbox{ for all } t\in [0,u]\}=1,
\qquad u>0.
 \eeqlb
For any $\phi \in C^2_\partial (\IR)$, we may use Lemma \ref{l4.6} to see that
 \beqnn
M^u_t(\phi)
:=
\<\hat P_{u-t}\phi, w_t\> - \<\hat P_t\hat P_{u-t}\phi, \mu\>
=
\int_0^t\int_{\hat\IR}\hat P_{u-s}\phi M(ds,dx),
\quad t\in [0,u],
 \eeqnn
is a continuous martingale with quadratic variation process
 \beqnn
\<M^u(\phi)\>_t
&=&
\int_0^t\<\sigma(\hat P_{u-s}\phi)^2, w_s\>ds
+ \int_0^t ds\int_{\hat\IR} \<h(z - \cdot)
\hat P_{u-s}(\phi^\prime), w_s\>^2 dz         \\
&=&
\int_0^t\<\sigma(\hat P_{u-s}\phi)^2, w_s\>ds
+ \int_0^t ds\int_{\hat\IR} \<h(z - \cdot)
(\hat P_{u-s}\phi)^\prime, w_s\>^2 dz.
 \eeqnn
By a martingale inequality we have
 \beqnn
& &\Q_\mu \bigg\{\sup_{0\le t\le u} |\<\hat P_{u-t}\phi, w_t\>
- \<\hat P_u \phi, \mu\>|^2\bigg\}   \\
&\le&
4\int_0^u\Q_\mu\{\<\sigma(\hat P_{u-s}\phi)^2, w_s\>\}ds
+ 4\int_0^u ds\int_{\hat\IR} \Q_\mu\{\<h(z - \cdot)
\hat P_{u-s}(\phi^\prime), w_s\>^2\} dz   \\
&\le&
4\int_0^u \<\sigma(\hat P_{u-s}\phi)^2, \mu\hat P_s\> ds
+ 4\int_{\hat\IR} h(z)^2 dz \int_0^u \Q_\mu\{\<1, w_s\>
\<[\hat P_{u-s}(\phi^\prime)]^2, w_s\>\}ds    \\
&\le&
4\int_0^u \<\sigma(\hat P_{u-s}\phi)^2, \mu\hat P_s\> ds
+ 4\|\phi^\prime\|^2 \int_{\hat\IR} h(z)^2 dz
\int_0^u \Q_\mu\{\<1, w_s\>^2\}ds.
 \eeqnn
Choose a sequence $\{\phi_k\} \subset C^2_\partial (\IR)$ such that $\phi_k (\cdot)
\to 1_{\{\partial\}}(\cdot)$ boundedly and $\|\phi_k^\prime\| \to 0$ as $k \to
\infty$. Replacing $\phi$ by $\phi_k$ in the above and letting $k \to \infty$ we
obtain (\ref{4.9}).
\qed

Combining Theorems \ref{t2.2} and \ref{t4.1} we get the existence of the SDSM in the
case where $\sigma \in C(\IR)^+$ extends continuously to $\hat \IR$.


\section{Measurable branching density}
\setcounter{equation}{0}
In this section, we shall use the dual process to extend the construction of the SDSM
to a general bounded Borel branching density. Given $\sigma\in B(\IR)^+$, let
$\{(M_t,Y_t): t\ge0\}$ be defined as in section 2. Choose any sequence of functions
$\{\sigma_k\} \subset C(\IR)^+$ which extends continuously to $\hat \IR$ and
$\sigma_k \to \sigma$ boundedly and pointwise. Suppose that $\{\mu_k\} \subset
M(\IR)$ and $\mu_k \to \mu\in M(\IR)$ as $k\to \infty$. For each $k\ge1$, let
$\{X_t^{(k)}: t\ge0\}$ be a SDSM with parameters $(a, \rho, \sigma_k)$ and initial
state $\mu_k \in M(\IR)$ and let $\Q_k$ denote the distribution of $\{X_t^{(k)}:
t\ge0\}$ on $C([0,\infty), M(\IR))$.

\blemma\label{l5.1}
Under the above hypotheses, $\{\Q_k\}$ is a tight sequence of probability measures on
$C([0,\infty), M(\IR))$.
\elemma

{\it Proof.} Since $\{\<1, X_t^{(k)}\>: t\ge0\}$ is a martingale, one can see as
in the proof of Lemma \ref{l4.1} that $\{X_t^{(k)}: t\ge0\}$ satisfies the compact
containment condition of Ethier and Kurtz (1986, p.142). Let $\L_k$ denote the
generator of $\{X_t^{(k)}: t\ge0\}$ and let $F$ be given by (\ref{4.1}) with $f\in
C^2_0 (\IR^n)$ and with $\{\phi_i\}\subset C^2_\partial (\IR)$. Then
 \beqnn
F(X_t^{(k)}) - F(X_0^{(k)}) - \int_0^t \L_k F(X_s^{(k)})ds,
\qquad t\ge0,
 \eeqnn
is a martingale. Since the sequence $\{\sigma_k\}$ is uniformly bounded, the tightness
of $\{X_t^{(k)}: t\ge0\}$ in $C([0,\infty), M(\hat\IR))$ follows from Lemma \ref{l4.4}
and the result of Ethier and Kurtz (1986, p.145). We shall prove that any limit point
of $\{\Q_k\}$ is supported by $C([0,\infty), M(\IR))$ so that $\{\Q_k\}$ is also tight
as probability measures on $C([0,\infty), M(\IR))$. Without loss of generality, we may
assume $\Q_k$ converges as $k\to \infty$ to $\Q_\mu$ by weak convergence of
probability measures on $C([0,\infty), M(\hat\IR))$. Let $\phi_n \in C^2(\IR)^+$ be such
that $\phi_n(x) =0$ when $\|x\|\le n$ and $\phi_n(x) =1$ when $\|x\|\ge 2n$ and
$\|\phi_n^\prime\| \to 0$ as $n\to \infty$. Fix $u>0$ and let $m_n$ be such that
$\phi_{m_n}(x) \le 2P_t \phi_n(x)$ for all $0\le t\le u$ and $x\in\IR$. For any
$\alpha>0$, the paths $w\in C([0,\infty), M(\hat\IR))$ satisfying $\sup_{0\le t\le u}
\<\phi_{m_n},w_t\> >\alpha$ constitute an open subset of $C([0,\infty), M(\hat\IR))$.
Then, by an equivalent condition for weak convergence,
 \beqnn
& & \Q_\mu\bigg\{\sup_{0\le t\le u} w_t(\{\partial\})>\alpha\bigg\}
\le
\Q_\mu\bigg\{\sup_{0\le t\le u} \<\phi_{m_n},w_t\> >\alpha\bigg\}    \\
&\le&
\liminf_{k\to \infty} \Q_k\bigg\{\sup_{0\le t\le u}
\<\phi_{m_n},w_t\> >\alpha\bigg\}
\le
\sup_{k\ge 1} \,\frac{4}{\alpha^2} \Q_k\bigg\{\sup_{0\le t\le u}
\<P_{u-t} \phi_{m_n},w_t\>^2\bigg\}    \\
&\le&
\sup_{k\ge 1} \,\frac{8}{\alpha^2} \Q_k\bigg\{\sup_{0\le t\le u}
|\<P_{u-t} \phi_{m_n},w_t\> - \<P_u\phi_{m_n},\mu_k\>|^2\bigg\}
+ \sup_{k\ge 1}\sup_{0\le t\le u} \,\frac{8}{\alpha^2}
\<P_u\phi_{m_n},\mu_k\>^2.
 \eeqnn
As in the proof of Theorem \ref{t4.1}, one can see that the right hand side goes to
zero as $m_n\to \infty$. Then $\Q_\mu$ is supported by $C([0,\infty), M(\IR))$.
\qed

\btheorem\label{t5.1}
The distribution $Q^{(k)}_t (\mu_k,\cdot)$ of $X_t^{(k)}$ on $M(\IR)$ converges as
$k\to \infty$ to a probability measure $Q_t(\mu,\cdot)$ on $M(\IR)$ given by
 \beqlb\label{5.1}
\int_{M(\IR)}\<f, \nu^m\>Q_t(\mu,d\nu)
=
\E^\sigma_{m,f} \bigg[\<Y_t, \mu^{M_t}\>
\exp\bigg\{\frac{1}{2}\int_0^t M_s(M_s-1)ds \bigg\}\bigg].
 \eeqlb
Moreover, $(Q_t)_{t\ge0}$ is a transition semigroup on $M(\IR)$.
\etheorem

{\it Proof.} By Lemma \ref{l5.1}, $\{Q^{(k)}_t (\mu_k,d\nu)\}$ is a tight sequence of
probability measures on $M(\IR)$. Take any subsequence $\{k_i\}$ so that $Q^{(k_i)}_t
(\mu_{k_i},d\nu)$ converges as $i\to\infty$ to some probability measure $Q_t(\mu,d\nu)$
on $M(\IR)$. By Lemma \ref{l2.1} we have
 \beqnn
& &\int_{M(\IR)}1_{[a,\infty)}(\<1,\nu\>) \<1, \nu^m\>Q^{(k)}_t(\mu_k,d\nu)   \\
&\le&
\frac{1}{\,a\,}\int_{M(\IR)} \<1, \nu^{m+1}\> Q^{(k)}_t(\mu_k,d\nu)  \\
&\le&
\frac{1}{\,a\,} \sum_{i=0}^m 2^{-i}(m+1)^i m^i \|\sigma_k\|^i \<1,\mu_k\>^{m-i+1},
 \eeqnn
which goes to zero as $a\to \infty$ uniformly in $k\ge1$. Then for $f\in C
(\hat\IR)^+$ we may regard $\{\<f,\nu^m\> Q^{(k)}_t (\mu_k,d\nu)\}$ as a tight sequence
of finite measures on $M(\hat\IR)$. By passing to a smaller subsequence $\{k_i\}$ we
may assume that $\<f,\nu^m\> Q^{(k_i)}_t (\mu_{k_i},d\nu)$ converges to a finite
measure $K_t(\mu,d\nu)$ on $M(\hat\IR)$. Then we must have $K_t(\mu,d\nu) = \<f,\nu^m\>
Q_t (\mu,d\nu)$. By Lemma \ref{l2.2} and the proof of Theorem \ref{l2.2}, $Q_t
(\mu,\cdot)$ is uniquely determined by (\ref{5.1}). Therefore, $Q^{(k)}_t
(\mu_k,\cdot)$ converges to $Q_t(\mu,\cdot)$ as $k\to \infty$. From the calculations
 \beqnn
& &\int_{M(\IR)}Q_r(\mu,d\eta) \int_{M(\IR)}\<f, \nu^m\>Q_t(\eta,d\nu) \\
&=&
\int_{M(\IR)}\E^\sigma_{m,f} \bigg[\<Y_t, \eta^{M_t}\>
\exp\bigg\{\frac{1}{2}\int_0^t M_s(M_s-1)ds \bigg\}\bigg] Q_r(\mu,d\eta)  \\
&=&
\E^\sigma_{m,f} \bigg[\int_{M(\IR)}\<Y_t, \eta^{M_t}\> Q_r(\mu,d\eta)
\exp\bigg\{\frac{1}{2}\int_0^t M_s(M_s-1)ds \bigg\}\bigg]  \\
&=&
\E^\sigma_{m,f} \bigg[\E^\sigma_{M_t,Y_t}
\bigg(\<Y_r, \mu^{M_r}\>
\exp\bigg\{\frac{1}{2}\int_0^r M_s(M_s-1)ds \bigg\}\bigg)
\exp\bigg\{\frac{1}{2}\int_0^t M_s(M_s-1)ds \bigg\}\bigg]  \\
&=&
\E^\sigma_{m,f} \bigg[\<Y_{r+t}, \mu^{M_{r+t}}\>
\exp\bigg\{\frac{1}{2}\int_0^{r+t} M_s(M_s-1)ds \bigg\}\bigg]  \\
&=&
\int_{M(\IR)}\<f, \nu^m\>Q_{r+t}(\eta,d\nu)
 \eeqnn
we have the Chapman-Kolmogorov equation. \qed

The existence of a SDSM with a general bounded measurable branching density function
$\sigma\in B(\IR)$ is given by the following

\btheorem\label{t5.2}
The sequence $\Q_k$ converges as $k\to \infty$ to a probability measure $\Q_\mu$ on
$C([0,\infty), M(\IR))$ under which the coordinate process $\{w_t: t\ge 0\}$ is a
diffusion with transition semigroup $(Q_t)_{t\ge0}$ defined by (\ref{5.1}). Let ${\cal
D}(\L)$ be the union of all functions of the form (\ref{4.1}) with $f\in C^2
(\IR^n)$ and $\{\phi_i\}\subset C^2 (\IR)$ and all functions of the form $F_{m,f}(\mu) =
\<f,\mu^m\>$ with $f\in C^2 (\IR^m)$. Then $\{w_t: t\ge0\}$ under $\Q_\mu$ solves the
$({\cal L}, {\cal D}(\L), \mu)$-martingale problem.
\etheorem

{\it Proof.} Let $\Q_\mu$ be the limit point of any subsequence $\{\Q_{k_i}\}$ of
$\{\Q_k\}$. Using Skorokhod's representation, we may construct processes
$\{X_t^{(k_i)}: t\ge0\}$ and $\{X_t: t\ge0\}$ with distributions $\Q_{k_i}$ and
$\Q_\mu$ on $C([0,\infty), M(\IR))$ such that $\{X_t^{(k_i)}: t\ge0\}$ converges
to $\{X_t: t\ge0\}$ a.s. when $i\to\infty$; see Ethier and Kurtz (1986, p.102). For
any $\{H_j\} _{j=1}^{n+1} \subset C(M(\hat\IR))$ and $0\le t_1 <\cdots <t_n <t_{n+1}$
we may use Theorem \ref{t5.1} and dominated convergence to see that
 \beqnn
& &\E\bigg\{\prod_{j=1}^n H_j(X_{t_j})H_{n+1}(X_{t_{n+1}})\bigg\}    \\
&=&
\lim_{i\to\infty} \E\bigg\{\prod_{j=1}^n H_j(X_{t_j}^{(k_i)})
H_{n+1}(X_{t_{n+1}}^{(k_i)})\bigg\}  \\
&=&
\lim_{i\to\infty} \E\bigg\{\prod_{j=1}^n H_j(X_{t_j}^{(k_i)}) \int_{M(\IR)}
H_{n+1}(\nu) Q_{t_{n+1}-t_n}^{(k_i)}(X_{t_n}^{(k_i)},d\nu)\bigg\}  \\
&=&
\E\bigg\{\prod_{j=1}^n H_j(X_{t_j}) \int_{M(\IR)}
H_{n+1}(\nu) Q_{t_{n+1}-t_n}(X_{t_n},d\nu)\bigg\}.
 \eeqnn
Then $\{X_t: t\ge0\}$ is a Markov process with transition semigroup $(Q_t)_{t\ge0}$
and actually $\Q_k \to \Q_\mu$ as $k\to \infty$. The strong Markov property holds
since $(Q_t)_{t\ge0}$ is Feller by (\ref{5.1}). To see the last assertion, one may
simply check that $({\cal L}, {\cal D}(\L))$ is a restriction of the generator of
$(Q_t)_{t\ge0}$. \qed


\section{Rescaled limits}
\setcounter{equation}0
In this section, we study the rescaled limits of the SDSM constructed in the last
section. Given any $\theta>0$, we defined the operator $K_\theta$ on $M(\IR)$ by
$K_\theta \mu(B) = \mu (\{\theta x: x\in B\})$. For a function $h\in B(\IR)$ we let
$h_\theta(x) = h(\theta x)$.

\blemma\label{l6.1}
Suppose that $\{X_t: t\ge0\}$ is a SDSM with parameters $(a, \rho, \sigma)$. Let
$X^\theta_t = \theta^{-2} K_\theta X_{\theta^2t}$. Then $\{X^\theta_t: t\ge0\}$ is a
SDSM with parameters $(a_\theta, \rho_\theta, \sigma_\theta)$.
\elemma

{\it Proof.} We shall compute the generator of $\{X^\theta_t: t\ge0\}$. Let $F(\mu)
= f(\<\phi,\mu\>)$ with $f\in C^2(\IR)$ and $\phi\in C^2(\IR)$. Note that
$F\circ K_\theta(\mu) = F(K_\theta\mu) = f(\<\phi_{1/\theta}, \mu\>)$. By the theory
of transformations of Markov processes, $\{K_\theta X_t: t\ge0\}$ has generator
$\L^\theta$ such that $\L^\theta F(\mu) = \L (F\circ K_\theta)(K_{1/\theta}\mu)$. Since
 \beqnn
\frac{d}{dx}\phi_{1/\theta}(x)
=
\frac{1}{\theta}(\phi^\prime)_{1/\theta}(x)
\quad\mbox{and}\quad
\frac{d^2}{dx^2}\phi_{1/\theta}(x)
=
\frac{1}{\theta^2}(\phi^{\prime\prime})_{1/\theta}(x),
 \eeqnn
it is easy to check that
 \beqnn
\L^\theta F(\mu)
&=&
\frac{1}{2\theta^2}f^\prime(\<\phi,\mu\>) \<a_{\theta}\phi^{\prime\prime},\mu\>  \\
& &\quad
+\frac{1}{2\theta^2} f^{\prime\prime}(\<\phi,\mu\>) \int_{\IR^2}\rho_{\theta}(x-y)
\phi^\prime(x) \phi^\prime(y) \mu(dx)\mu(dy)    \\
& &\quad
+\frac{1}{2} f^{\prime\prime}(\<\phi,\mu\>)\<\sigma_{\theta} \phi^2,\mu\>.
 \eeqnn
Then one may see that $\{\theta^{-2}K_\theta X_t: t\ge0\}$ has generator $\L_\theta$
such that
 \beqnn
\L_\theta F(\mu)
&=&
\frac{1}{2\theta^2}f^\prime(\<\phi,\mu\>) \<a_{\theta}\phi^{\prime\prime},\mu\> \\
& &\quad
+\frac{1}{2\theta^2} f^{\prime\prime}(\<\phi,\mu\>) \int_{\IR^2}\rho_{\theta}(x-y)
\phi^\prime(x) \phi^\prime(y) \mu(dx)\mu(dy)    \\
& &\quad
+\frac{1}{2\theta^2} f^{\prime\prime}(\<\phi,\mu\>)\<\sigma_{\theta} \phi^2,\mu\>,
 \eeqnn
and hence $\{X^\theta_t: t\ge0\}$ has the right generator $\theta^2 \L_\theta$. \qed

\btheorem\label{t6.1}
Suppose that $(\itOmega, X_t, \Q_\mu)$ is a realization of the SDSM with parameters
$(a,\rho,\sigma)$ with $|c(x)| \ge \epsilon >0$ for all $x\in \IR$. Then there is a
$\lambda \times \lambda \times \Q_\mu$-measurable function $X_t(\omega,x)$ such that
$\Q_\mu\{\omega\in \itOmega: X_t(\omega,dx)$ is absolutely continuous with respect to
the Lebesgue measure with density $X_t(\omega,x)$ for $\lambda$-a.e. $t>0\} =1$.
Moreover, for $\lambda \times \lambda$-a.e. $(t,x) \in [0,\infty) \times \IR$ we have
 \beqlb\label{6.1}
\Q_\mu\{X_t(x)^2\}
&=&
\int_{\IR^2}p_t^2(y,z;x,x)\mu(dx)\mu(dy) \nonumber  \\
& &\quad + \int_0^tds \int_{\IR}\mu(dy) \int_{\IR}\sigma(z)p_s^2(z,z;x,x)
p_{t-s}(y,z)dz.
 \eeqlb
\etheorem

{\it Proof.} Recall (\ref{1.9}). For $r_1>0$ and $r_2>0$ we use (\ref{2.7}) and
(\ref{5.1}) to see that
 \beqnn
& &\Q_\mu\{\<g^1_{\epsilon r_1}(x,\cdot),X_t\>
\<g^1_{\epsilon r_2}(x,\cdot),X_t\>\}
=
\Q_\mu\{\<g^1_{\epsilon r_1}(x,\cdot)\otimes
g^1_{\epsilon r_2}(x,\cdot), X_t^2\>\}  \\
&=&
\<P_t^2 g^1_{\epsilon r_1}(x,\cdot)\otimes
g^1_{\epsilon r_2}(x,\cdot), \mu^2\>
+ \int_0^t\<P_{t-s}\itPhi_{12}P_s^2g^1_{\epsilon r_1}
(x,\cdot)\otimes g^1_{\epsilon r_2}(x,\cdot),\mu\>ds  \\
&=&
\int_{\IR^2}P_t^2 g^1_{\epsilon r_1}(x,\cdot) \otimes
g^1_{\epsilon r_2}(x,\cdot)(y,z)\mu(dy)\mu(dz)  \\
& & + \int_0^tds\int_{\IR}\mu(dy)\int_{\IR} \sigma(z)
P_s^2 g^1_{\epsilon r_1}(x,\cdot) \otimes  g^1_{\epsilon r_2}
(x,\cdot)(z,z)p_{t-s}(y,z)dz.
 \eeqnn
Observe that
 \beqnn
P_t^2 g^1_{\epsilon r_1}(x,\cdot) \otimes g^1_{\epsilon r_2}(x,\cdot)(y,z)
=
\int_{\IR^2}g^1_{\epsilon r_1}(x,z_1)g^1_{\epsilon r_2}(x,z_2)
p_t^2(y,z;z_1,z_2)dz_1dz_2
 \eeqnn
converges to $p_t^2(y,z;x,x)$ boundedly as $r_1\to 0$ and $r_2\to 0$. Note also that
 \beqnn
& &\int_{\IR} \sigma(z)P_s^2 g^1_{\epsilon r_1}(x,\cdot) \otimes
g^1_{\epsilon r_2}(x,\cdot)(z,z) p_{t-s}(y,z) dz  \\
&\le&
\mbox{const} \cdot \|\sigma\|\frac{1}{\sqrt{s}}\int_{\IR} T_{\epsilon s}
g^1_{\epsilon r_1}(x;\cdot)(z) g^1_{\epsilon (t-s)}(y,z) dz   \\
&\le&
\mbox{const} \cdot \|\sigma\|\frac{1}{\sqrt{s}}g^1_{\epsilon (t+r_1)}(y,x)   \\
&\le&
\mbox{const} \cdot \|\sigma\|\frac{1}{\sqrt{st}}.
 \eeqnn
By dominated convergence theorem we get
 \beqnn
& &\lim_{r_1,r_2\to 0}\Q_\mu\{\<g^1_{\epsilon r_1}(x,\cdot),X_t\>
\<g^1_{\epsilon r_2}(x,\cdot),X_t\>\}  \\
&=&
\int_{\IR^2} p_t^2(y,z;x,x)\mu(dy)\mu(dz)  \\
& &\qquad + \int_0^tds\int_{\IR}\mu(dy)\int_{\IR} \sigma(z) p_t^2(z,z;x,x)
p_{t-s}(y,z) dz.
 \eeqnn
Then it is easy to check that
 \beqnn
\lim_{r_1,r_2\to0}\int_0^Tdt\int_{\IR} \Q_\mu\{\<g^1_{\epsilon r_1}(x,\cdot)
- g^1_{\epsilon r_2}(x,\cdot),X_t\>^2\}dx =0
 \eeqnn
for each $T>0$, so there is a $\lambda \times \lambda \times \Q_\mu$-measurable
function $X_t(\omega,x)$ satisfying (\ref{6.1}) and
 \beqlb\label{6.2}
\lim_{r\to0} \int_{\IR} g^1_{\epsilon r}(x,y) X_t(\omega,dy)
=
X_t(\omega,x)
 \eeqlb
in $L^2(\lambda \times \lambda \times \Q_\mu)$. For any square integrable $\phi\in
C(\IR)$,
 \beqlb\label{6.3}
& &\int_0^T \Q_\mu\bigg\{\bigg|\<\phi,X_t\>
- \int_{\IR}\phi(x)X_t(x)dx\bigg|^2\bigg\}dt \nonumber \\
&\le&
2 \int_0^T \Q_\mu\left\{\<\phi-T_{\epsilon r}\phi,X_t\>^2\right\}dt  \nonumber \\
& &\qquad
+ 2\int_0^T \Q_\mu\bigg\{\bigg|\<T_{\epsilon r}\phi,X_t\>
- \int_{\IR}\phi(x)X_t(x)dx\bigg|^2\bigg\}dt,
 \eeqlb
and by Schwarz inequality,
 \beqnn
& &\Q_\mu\bigg\{\bigg|\<T_{\epsilon r}\phi,X_t\>
- \int_{\IR}\phi(x)X_t(x)dx\bigg|^2\bigg\}    \\
&=&
\Q_\mu\bigg\{\bigg|\int_{\IR}X_t(dx)\int_{\IR} \phi(x)g^1_{\epsilon r}(y,x) dx
- \int_{\IR}\phi(x)X_t(x)dx\bigg|^2\bigg\}   \\
&=&
\Q_\mu\bigg\{\bigg|\int_{\IR}[\<g^1_{\epsilon r}(\cdot,x),X_t\>
- X_t(x)]\phi(x)dx\bigg|^2\bigg\}   \\
&\le&
\int_{\IR}\Q_\mu\left\{|\<g^1_{\epsilon r}(\cdot,x),X_t\> - X_t(x)|^2\right\}dx
\int_{\IR}\phi(x)^2 dx.
 \eeqnn
By this and (\ref{6.2}) we get
 \beqnn
\lim_{r\to 0}\int_0^T \Q_\mu\bigg\{\bigg|\<T_{\epsilon r}\phi,X_t\>
- \int_{\IR}\phi(x)X_t(x)dx\bigg|^2\bigg\}dt
=0.
 \eeqnn
On the other hand, using (\ref{2.8}) and (\ref{5.1}) one may see that
 \beqnn
\lim_{r\to 0}\Q_\mu \{\<\phi - T_{\epsilon r}\phi,X_t\>^2\}
\le
\lim_{r\to 0} \|\phi - T_{\epsilon r}\phi\|^2 \Q_\mu\{\<1,X_t\>^2\}
= 0.
 \eeqnn
Then letting $r\to0$ in (\ref{6.3}) we have
 \beqnn
\int_0^T \Q_\mu\bigg\{\bigg|\<\phi,X_t\>
- \int_{\IR}\phi(x)X_t(x)dx\bigg|^2\bigg\}dt
=0,
 \eeqnn
completing the proof.
\qed

By Theorem \ref{t6.1}, for $\lambda \times \lambda$-a.e. $(t,x) \in [0,\infty)\times
\IR$ we have
 \beqlb\label{6.4}
\Q_\mu\{X_t(x)^2\}
&\le&
\mbox{const}\cdot\bigg[\frac{1}{\sqrt{t}}\<1,\mu\>\int_{\IR}g^1_{\epsilon
t}(x,y)\mu(dy)  \nonumber \\
& &\qquad
+ \int_0^t\frac{ds}{\sqrt{s}}\int_{\IR}\mu(dy) \int_{\IR}\|\sigma\|
g^1_{\epsilon s}(z,x) g^1_{\epsilon (t-s)}(z,x)dz\bigg] \nonumber \\
&\le&
\mbox{const}\cdot\bigg[\frac{1}{\sqrt{t}}\<1,\mu\> + \sqrt{t}\|\sigma\|\bigg]
\int_{\IR}g^1_{\epsilon t}(x,y)\mu(dy).
 \eeqlb

\btheorem\label{t6.2}
Suppose that $\{X_t: t\ge0\}$ is a SDSM with parameters $(a, \rho, \sigma)$ with
$|c(x)| \ge \epsilon >0$ for all $x\in \IR$. Let $X^\theta_t = \theta^{-2} K_\theta
X_{\theta^2t}$. Assume $a(x) \to a_\partial$, $\sigma(x) \to \sigma_\partial$ and
$\rho(x) \to 0$ as $|x| \to \infty$. Then the conditional distribution of
$\{X^\theta_t: t\ge0\}$ given $X^\theta_0 = \mu\in M(\IR)$ converges as $\theta \to
\infty$ to that of a super Brownian motion with underlying generator $(a_\partial/2)
\itDelta$ and uniform branching density $\sigma_\partial$.
\etheorem

{\it Proof.} Since $\|\sigma_\theta\| = \|\sigma\|$ and $X^\theta_0 = \mu$, as in the
proof of Lemma \ref{l5.1} one can see that the family $\{X^\theta_t: t\ge0\}$ is tight in
$C([0,\infty), M(\IR))$. Choose any sequence $\theta_k\to \infty$ such that the
distribution of $\{X^{\theta_k}_t: t\ge0\}$ converges to some probability measure
$\Q_\mu$ on $C([0,\infty), M(\IR))$. We shall prove that $\Q_\mu$ is the solution of
the martingale problem for the super Brownian motion so that actually the distribution
of $\{X^{\theta}_t: t\ge0\}$ converges to $\Q_\mu$ as $\theta \to \infty$. By
Skorokhod's representation, we can construct processes $\{X^{(k)}_t: t\ge0\}$ and
$\{X^{(0)}_t: t\ge0\}$ such that $\{X^{(k)}_t: t\ge0\}$ and $\{X^{\theta_k}_t:
t\ge0\}$ have identical distributions, $\{X^{(0)}_t: t\ge0\}$ has the distribution
$\Q_\mu$ and $\{X^{(k)}_t: t\ge0\}$ converges a.s. to $\{X^{(0)}_t: t\ge0\}$ in
$C([0,\infty), M(\IR))$. Let $F(\mu) = f(\<\phi,\mu\>)$ with $f\in C^2(\IR)$ and
$\phi \in C^2 (\IR)$. Then for each $k\ge0$,
 \beqlb\label{6.5}
F(X_t^{(k)}) - F(X_0^{(k)}) - \int_0^t \L_k F(X_s^{(k)})ds,
\qquad t\ge0,
 \eeqlb
is a martingale, where $\L_k$ is given by
 \beqnn
\L_k F(\mu)
&=&
\frac{1}{2}f^\prime(\<\phi,\mu\>) \<a_{\theta_k}\phi^{\prime\prime},\mu\>
+ \frac{1}{2} f^{\prime\prime}(\<\phi,\mu\>)\<\sigma_{\theta_k} \phi^2,\mu\>  \\
& & + \frac{1}{2} f^{\prime\prime}(\<\phi,\mu\>) \int_{\IR^2}\rho_{\theta_k}(x-y)
\phi^\prime(x) \phi^\prime(y) \mu(dx)\mu(dy).
 \eeqnn
Observe that
 \beqnn
& &\int_0^t\E\{|f^{\prime}(\<\phi,X_s^{(k)}\>)|
\<|a_{\theta_k} - a_\partial|\phi^{\prime\prime},X_s^{(k)}\>\}ds  \\
&\le&
\|f^{\prime}\|\|\phi^{\prime\prime}\|\int_0^t
\E\{\<|a_{\theta_k} - a_\partial|, X_s^{(k)}\>\}ds  \\
&\le&
\|f^{\prime}\|\|\phi^{\prime\prime}\|\int_0^t
\<P_s|a_{\theta_k} - a_\partial|, \mu\> ds  \\
&\le&
\|f^{\prime}\|\|\phi^{\prime\prime}\|\int_0^tds \int_{\IR} \mu(dx) \int_{\IR}
|a_{\theta_k}(y) - a_\partial| p_s(x,y)dy.
 \eeqnn
Then we have
 \beqlb\label{6.6}
\lim_{k\to\infty} \int_0^t\E\{|f^{\prime}(\<\phi,X_s^{(k)}\>)|
\<|a_{\theta_k} - a_\partial|\phi^{\prime\prime},X_s^{(k)}\>\}ds
=0.
 \eeqlb
In the same way, one sees that
 \beqlb\label{6.7}
\lim_{k\to\infty} \int_0^t\E\{|f^{\prime\prime}(\<\phi,X_s^{(k)}\>)|
\<|\sigma_{\theta_k} - \sigma_\partial|\phi^2,X_s^{(k)}\>\}ds
=0.
 \eeqlb
Using the density process of $\{X^{(k)}_t: t\ge0\}$ we have the following estimates
 \beqnn
& &\E\bigg|f^{\prime\prime}(\<\phi,X_s^{(k)}\>) \int_{\IR^2}\rho_{\theta_k}(x-y)
\phi^\prime(x) \phi^\prime(y) X_s^{(k)}(dx)X_s^{(k)}(dy)\bigg|     \\
&\le&
\|f^{\prime\prime}\|\int_{\IR^2}|\rho_{\theta_k}(x-y)| |\phi^\prime(x)\phi^\prime(y)|
\E\{X_s^{(k)}(x)X_s^{(k)}(y)\}dxdy    \\
&\le&
\|f^{\prime\prime}\|\int_{\IR^2}|\rho_{\theta_k}(x-y)| |\phi^\prime(x)\phi^\prime(y)|
\E\{X_s^{(k)}(x)^2\}^{1/2}\E\{X_s^{(k)}(y)^2\}^{1/2}dxdy    \\
&\le&
\|f^{\prime\prime}\|\bigg(\int_{\IR^2}|\rho_{\theta_k}(x-y)|^2
|\phi^\prime(x)\phi^\prime(y)|^2 dxdy
\int_{\IR^2}\E\{X_s^{(k)}(x)^2\}\E\{X_s^{(k)}(y)^2\}dxdy\bigg)^{1/2} \\
&\le&
\|f^{\prime\prime}\|\bigg(\int_{\IR^2}|\rho_{\theta_k}(x-y)|^2
|\phi^\prime(x)\phi^\prime(y)|^2 dxdy\bigg)^{1/2}
\int_{\IR}\E\{X_s^{(k)}(x)^2\}dx.
 \eeqnn
By (\ref{6.4}), for any fixed $t\ge0$,
 \beqnn
\int_0^tds\int_{\IR}\E\{X_s^{(k)}(x)^2\}dx
 \eeqnn
is uniformly bounded in $k\ge1$. Since $\rho_{\theta_k}(x-y) \to 0$ for $\lambda
\times \lambda$-a.e. $(x,y)\in \IR^2$ and since $\|\rho_{\theta_k}\| = \|\rho\|$, we
have
 \beqnn
\lim_{k\to\infty}\int_{\IR^2} |\rho_{\theta_k}(x-y)|^2
|\phi^\prime(x)\phi^\prime(y)|^2 dxdy
=0
 \eeqnn
when $\phi^\prime\in L^2(\lambda)$. Then
 \beqlb\label{6.8}
\lim_{k\to\infty}\E\bigg|f^{\prime\prime}(\<\phi,X_s^{(k)}\>) \int_{\IR^2}
\rho_{\theta_k}(x-y) \phi^\prime(x) \phi^\prime(y) X_s^{(k)}(dx)X_s^{(k)}(dy)\bigg|
=0.
 \eeqlb
Using (\ref{6.6}),(\ref{6.7}), (\ref{6.8}) and the martingale property of (\ref{6.5})
ones sees in a similar way as in the proof of Lemma \ref{l4.2} that
 \beqnn
F(X_t^{(0)}) - F(X_0^{(0)}) - \int_0^t \L_0 F(X_s^{(0)})ds,
\qquad t\ge0,
 \eeqnn
is a martingale, where $\L_0$ is given by
 \beqnn
\L_0 F(\mu)
=
\frac{1}{2}a_\partial f^\prime(\<\phi,\mu\>) \<\phi^{\prime\prime},\mu\>
+ \frac{1}{2}\sigma_\partial f^{\prime\prime}(\<\phi,\mu\>)\<\phi^2,\mu\>.
 \eeqnn
This clearly implies that $\{X^{(0)}_t: t\ge0\}$ is a solution of the martingale
problem of the super Brownian motion. \qed


\section{Measure-valued catalysts}
\setcounter{equation}0

In this section, we assume $|c(x)| \ge \epsilon>0$ for all $x\in \IR$ and give
construction for a class of SDSM with measure-valued catalysts. We start from the
construction of a class of measure-valued dual processes. Let $M_B(\IR)$ denote the
space of Radon measures $\zeta$ on $\IR$ to which there correspond constants
$b(\zeta)>0$ and $l(\zeta)>0$ such that
 \beqlb\label{7.1}
\zeta([x,x+l(\zeta)]) \le b(\zeta)l(\zeta),
\qquad x\in \IR.
 \eeqlb
Clearly, $M_B(\IR)$ contains all finite measures and all Radon measures which are
absolutely continuous with respect to the Lebesgue measure with bounded densities.
Let $M_B(\IR^m)$ denote the space of Radon measures $\nu$ on $\IR^m$ such that
 \beqlb\label{7.2}
\nu(dx_1,\cdots,dx_m)
=
f(x_1,\cdots,x_m) dx_1,\cdots,dx_{m-1} \zeta(dx_m)
 \eeqlb
for some $f\in C(\IR^m)$ and $\zeta\in M_B(\IR)$. We endow $M_B(\IR^m)$ with the
topology of vague convergence. Let $M_A(\IR^m)$ denote the subspace of $M_B(\IR^m)$
comprising of measures which are absolutely continuous with respect to the Lebesgue
measure and have bounded densities. For $f\in C(\IR^m)$, we define $\lambda^m_f \in
M_A(\IR^m)$ by $\lambda^m_f(dx) = f(x)dx$. Let $\bM$ be the topological union of
$\{M_B(\IR^m): m= 1,2,\cdots\}$.

\blemma\label{l7.1}
If $\zeta \in M_B(\IR)$ satisfies (\ref{7.1}), then
 \beqnn
\int_{\IR} p_t(x,y) \zeta(dy)
\le
h(\epsilon,\zeta; t)/\sqrt{t},
\qquad t>0, x\in \IR,
 \eeqnn
where
 \beqnn
h(\epsilon,\zeta; t)
=
\mbox{const} \cdot b(\zeta)
\bigg[2l(\zeta) + \sqrt{2\pi \epsilon t}\bigg],
\qquad t>0.
 \eeqnn
\elemma

{\it Proof.} Using (\ref{1.9}) and (\ref{7.1}) we have
 \beqnn
\int_{\IR} p_t(x,y) \zeta(dy)
&\le&
\mbox{const}\cdot\int_{\IR} g_{\epsilon t}(x,y) \zeta (dy)   \\
&\le&
\mbox{const}\cdot\frac{2b(\zeta)l(\zeta)}{\sqrt{2\pi \epsilon t}}\sum_{k=0}^\infty
\exp\bigg\{-\frac{k^2l(\zeta)^2}{2\epsilon t}\bigg\}  \\
&\le&
\mbox{const}\cdot\frac{b(\zeta)}{\sqrt{2\pi \epsilon t}}
\bigg[2l(\zeta) + \int_{\IR} \exp\bigg\{-\frac{y^2}
{2\epsilon t}\bigg\} dy\bigg]  \\
&\le&
\mbox{const}\cdot \frac{b(\zeta)}{\sqrt{2\pi \epsilon t}}
\bigg[2l(\zeta) + \sqrt{2\pi \epsilon t}\bigg],
 \eeqnn
giving the desired inequality. \qed

Fix $\eta \in M_B(\IR)$ and let $\itPhi_{ij}$ be the mapping from $M_A(\IR^m)$
to $M_B(\IR^{m-1})$ defined by
 \beqlb\label{7.3}
& &\itPhi_{ij}\mu(dx_1,\cdots,dx_{m-1}) \nonumber \\
&=&
\mu^\prime(x_1,\cdots,x_{m-1},\cdots,x_{m-1},\cdots,x_{m-2})
dx_1 \cdots dx_{m-2} \eta(dx_{m-1}),
 \eeqlb
where $\mu^\prime$ denotes the Radon-Nikodym derivative of $\mu$ with respect to the
$m$-dimensional Lebesgue measure, and $x_{m-1}$ is in the places of the $i$th and the
$j$th variables of $\mu^\prime$ on the right hand side. We may also regard
$(P_t^m)_{t>0}$ as operators on $M_B(\IR^m)$ determined by
 \beqlb\label{7.4}
P^m_t\nu(dx) = \int_{\IR^m} p_t^m(x,y)\nu(dy)dx,
\qquad t>0, x\in \IR^m.
 \eeqlb
By Lemma \ref{l7.1} one can show that each $P_t^m$ maps $M_B(\IR^m)$ to $M_A(\IR^m)$
and, for $f\in C(\IR^m)$,
 \beqlb\label{7.5}
P^m_t\lambda^m_f(dx) = P^m_tf(x) dx,
\qquad t>0, x\in \IR^m.
 \eeqlb
Let $\{M_t: t\ge 0\}$ and $\{\itGamma_k: 1\le k\le M_0-1\}$ be defined as in section 2.
Then
 \beqlb\label{7.6}
Z_t
=
P^{M_{\tau_k}}_{t-\tau_k} \itGamma_k P^{M_{\tau_{k-1}}}_{\tau_k -\tau_{k-1}}
\itGamma_{k-1} \cdots P^{M_{\tau_1}}_{\tau_2 -\tau_1} \itGamma_1
P^{M_0}_{\tau_1}Z_0,
\quad \tau_k \le t < \tau_{k+1}, 0\le k\le M_0-1,
 \eeqlb
defines a Markov process $\{Z_t: t\ge0\}$ taking values from $\bM$. Of course,
$\{(M_t, Z_t): t\ge 0\}$ is also a Markov process. We shall suppress the dependence
of $\{Z_t: t\ge 0\}$ on $\eta$ and let $\E^\eta_{m,\nu}$ denote the expectation
given $M_0=m$ and $Z_0=\nu \in M_B(\IR^m)$. Observe that by (\ref{7.4}) and
(\ref{7.6}) we have
 \beqlb\label{7.7}
& &\E^{\eta}_{m,\nu} \bigg[\<Z_t^\prime, \mu^{M_t}\>
\exp\bigg\{\frac{1}{2}\int_0^t M_s(M_s-1)ds \bigg\}\bigg]  \nonumber \\
&=&
\<(P^m_t\nu)^\prime,\mu^m\> \\
&+&
\frac{1}{2}\sum_{i,j=1, i\neq j}^m \int_0^t
\E^{\eta}_{m-1,\itPhi_{ij}P^m_u\nu} \bigg[\<Z_{t-u}^\prime,\mu^{M_{t-u}}\>
\exp\bigg\{\frac{1}{2}\int_0^{t-u} M_s(M_s-1)ds \bigg\}\bigg]du. \nonumber
 \eeqlb

\blemma\label{l7.2}
Let $\eta \in M_B(\IR)$. For any integer $k\ge1$, define $\eta_k \in M_A(\IR)$ by
 \beqnn
\eta_k(dx)
=
kl(\eta)^{-1}\eta ((il(\eta)/k, (i+1)l(\eta)/k])dx,
\qquad x\in (il(\eta)/k, (i+1)l(\eta)/k],
 \eeqnn
where $i=\cdots,-2,-1,0,1,2,\cdots$. Then $\eta_k \to \eta$ by weak convergence as
$k\to \infty$ and
 \beqnn
\eta_k([x,x+l(\eta)]) \le 2b(\eta)l(\eta),
\qquad x\in \IR.
 \eeqnn
\elemma

{\it Proof.} The convergence $\eta_k \to \eta$ as $k\to \infty$ is clear. For any
$x\in \IR$ there is an integer $i$ such that
 \beqnn
[x,x+l(\eta)]
\subset
(il(\eta)/k, (i+1)l(\eta)/k+l(\eta)].
 \eeqnn
Therefore, we have
 \beqnn
\eta_k([x,x+l(\eta)])
&\le&
\eta_k ((il(\eta)/k, (i+1)l(\eta)/k+l(\eta)])   \\
&=&
\eta ((il(\eta)/k, (i+1)l(\eta)/k+l(\eta)])   \\
&\le&
\eta ((il(\eta)/k, il(\eta)/k+2l(\eta)])   \\
&\le&
2b(\eta)l(\eta),
 \eeqnn
as desired.
\qed

\blemma\label{l7.3}
If $\eta \in M_B(\IR)$ and if $\nu \in M_B(\IR^m)$ is given by (\ref{7.2}), then
 \beqlb\label{7.8}
& &\E^{\eta}_{m,\nu} \bigg[\<Z_t^\prime, \mu^{M_t}\>
\exp\bigg\{\frac{1}{2}\int_0^t M_s(M_s-1)ds \bigg\}\bigg]   \nonumber \\
&\le&
\|f\| h(\epsilon,\zeta;t)\bigg[\<1,\mu\>^m/\sqrt{t}
+ \sum_{k=1}^{m-1}2^k m^k(m-1)^k
h(\epsilon,\eta;t)^k\<1,\mu\>^{m-k}t^{k/2}\bigg].
 \eeqlb
(Note that the left hand side of (\ref{7.8}) is well defined since $Z_t\in M_A(\IR)$
a.s. for each $t>0$ by (\ref{7.6}).)
\elemma

{\it Proof.} The left hand side of (\ref{7.8}) can be decomposed as $\sum_{k=0}^{m-1}
A_k$ with
 \beqnn
A_k = \E^{\eta}_{m,\nu} \bigg[\<Z_t^\prime, \mu^{M_t}\>
\exp\bigg\{\frac{1}{2}\int_0^t M_s(M_s-1)ds \bigg\}
1_{\{\tau_k\le t<\tau_{k+1}\}}\bigg].
 \eeqnn
By (\ref{7.2}) and Lemma \ref{l7.1},
 \beqnn
A_0
=
\<(P^m_t\nu)^\prime,\mu^m\>
\le
\|f\| h(\epsilon,\zeta;t)\<1,\mu\>^m/\sqrt{t}.
 \eeqnn
By the construction (\ref{7.6}) we have
 \beqnn
A_k
&=&
\frac{m!(m-1)!}{2^k (m-k)!(m-k-1)!}
\int_0^t ds_1 \int_{s_1}^t ds_2 \cdots
\int_{s_{k-1}}^t    \\
& &\qquad
\E^{\eta}_{m,\nu}\{\<(P^{m-k}_{t-s_k} \itGamma_k \cdots
P^{m-1}_{s_2 -s_1}\itGamma_1 P^m_{s_1}\nu)^\prime,\mu^{m-k}\>
|\tau_j = s_j: 1\le j\le k\} ds_k
 \eeqnn
for $1\le k\le m-1$. Observe that
 \beqlb\label{7.9}
\int_{s_{k-1}}^t\frac{ds_k}{\sqrt{t-s_k}\sqrt{s_k-s_{k-1}}}
\le
\frac{2\sqrt{2}}{\sqrt{t-s_{k-1}}}
\int_{(t+s_{k-1})/2}^t\frac{ds_k}{\sqrt{t-s_k}}
\le
\frac{4\sqrt{t}}{\sqrt{t-s_{k-1}}}.
 \eeqlb
By (\ref{7.5}) we have $P^{m-k}_s\lambda^{m-k}_h \le \lambda^{m-k}_{\|h\|}$ for
$h\in C(\IR^{m-k})$. Then using (\ref{7.9}) and Lemma \ref{l7.1} inductively we get
 \beqnn
A_k
&\le&
\frac{m!(m-1)!\|f\|}{2^k (m-k)!(m-k-1)!}\int_0^tds_1\int_{s_1}^tds_2
\cdots \int_{s_{k-1}}^t    \\
& &\qquad
\frac{h(\epsilon,\zeta;t) h(\epsilon,\eta;t)^k\<1,\mu\>^{m-k}}
{\sqrt{t-s_k} \cdots \sqrt{s_2-s_1}\sqrt{s_1}}ds_k  \\
&\le&
\frac{2^k m!(m-1)!\|f\|}{(m-k)!(m-k-1)!} h(\epsilon,\zeta;t)
h(\epsilon,\eta;t)^k\<1,\mu\>^{m-k}t^{k/2}  \\
&\le&
2^k m^k(m-1)^k \|f\| h(\epsilon,\zeta;t)
h(\epsilon,\eta;t)^k\<1,\mu\>^{m-k}t^{k/2}.
 \eeqnn
Returning to the decomposition we get the desired estimate. \qed

\blemma\label{l7.4}
Suppose $\eta \in M_B(\IR)$ and define $\eta_k \in M_A(\IR)$ as in Lemma \ref{l7.2}.
Assume that $\mu_k \to \mu$ weakly as $k\to \infty$. Then we have
 \beqnn
& &\E^\eta_{m,\nu} \bigg[\<Z_t^\prime, \mu^{M_t}\>
\exp\bigg\{\frac{1}{2}\int_0^t M_s(M_s-1)ds \bigg\}\bigg]  \nonumber  \\
&=&
\lim_{k\to \infty} \E^{\eta_k}_{m,\nu} \bigg[\<Z_t^\prime, \mu_k^{M_t}\>
\exp\bigg\{\frac{1}{2}\int_0^t M_s(M_s-1)ds \bigg\}\bigg].
 \eeqnn
\elemma

{\it Proof.} Based on (\ref{7.7}), the desired result follows by a similar argument
as in the proof of Lemma \ref{l2.2}.
\qed

Let $\eta\in M_B(\IR)$ and let $\eta_k$ be defined as in Lemma \ref{l7.2}. Let
$\sigma_k$ denote the density of $\eta_k$ with respect to the Lebesgue measure and let
$\{X_t^{(k)}: t\ge0\}$ be a SDSM with parameters $(a, \rho, \sigma_k)$ and initial
state $\mu_k\in M(\IR)$. Assume that $\mu_k \to \mu$ weakly as $k\to \infty$. Then we
have the following

\btheorem\label{t7.1}
The distribution $Q^{(k)}_t (\mu_k,\cdot)$ of $X_t^{(k)}$ on $M(\IR)$ converges as
$k\to \infty$ to a probability measure $Q_t(\mu,\cdot)$ on $M(\IR)$ given by
 \beqlb\label{7.10}
\int_{M(\IR)}\<f, \nu^m\>Q_t(\mu,d\nu)
=
\E^\eta_{m,\lambda^m_f} \bigg[\<Z_t^\prime, \mu^{M_t}\>
\exp\bigg\{\frac{1}{2}\int_0^t M_s(M_s-1)ds \bigg\}\bigg].
 \eeqlb
Moreover, $(Q_t)_{t\ge0}$ is a transition semigroup on $M(\IR)$.
\etheorem

{\it Proof.} With Lemmas \ref{l7.3} and \ref{l7.4}, this is similar to the proof of
Theorem \ref{t5.1}.
\qed

A Markov process with transition semigroup defined by (\ref{7.10}) is the so-called
SDSM with measure-valued catalysts.

\end{document}